\documentclass[a4paper,10pt]{article}

\usepackage{geometry}
\usepackage[english]{babel}
\usepackage{moreverb}

\usepackage{amsmath, bm}
\usepackage{caption}
\usepackage{booktabs}

\usepackage{hyperref}
\usepackage{color}
\usepackage{multirow}
\usepackage{mathrsfs} 
\usepackage{amssymb}
\usepackage{latexsym}
\usepackage{fancyhdr}
\usepackage{graphicx}
\usepackage{newlfont}
\usepackage{textcomp}
\usepackage{float}
\usepackage{tikz}
\usepackage{mleftright}
\usepackage{multirow}

\usepackage{enumitem}

\newcommand\BibTeX{{\rmfamily B\kern-.05em \textsc{i\kern-.025em b}\kern-.08em
T\kern-.1667em\lower.7ex\hbox{E}\kern-.125emX}}

\usepackage{etoolbox}
\AtBeginEnvironment{bmatrix}{\setlength{\arraycolsep}{1pt}} 
\makeatletter
\renewcommand*\env@matrix[1][*\c@MaxMatrixCols c]{%
  \hskip -\arraycolsep
  \let\@ifnextchar\new@ifnextchar
  \array{#1}}
\makeatother

\newdimen\nodeDist
\tikzset{
    position/.style args={#1:#2 from #3}{
        at=(#3.#1), anchor=#1+180, shift=(#1:#2)
    }
}

\newcommand{\new}[1]{\textcolor{black}{#1}}
\newcommand{\aled}[1]{\widehat{#1}}

\title{{\bf{A field-split preconditioning technique for fluid-structure interaction problems with applications in biomechanics}}}
\author{%
  S. Calandrini\footnote{Department of Scientific Computing, Florida State University, Tallahassee FL 32306, USA.}%
  \and E. Aulisa\footnote{Department of Mathematics and Statistics, Texas Tech University, Lubbock TX 79409, USA.}%
  \and G. Ke \footnote{Department of Mathematics and Physical Sciences, Louisiana State University at Alexandria, Alexandria LA 71302, USA.}%
  }
\date{}

\begin{document}
 
\maketitle

\begin{abstract} 
We present a novel preconditioning technique for Krylov subspace 
algorithms to solve fluid-structure interaction (FSI) linearized systems arising from finite element discretizations.
An outer Krylov subspace solver preconditioned with a geometric multigrid (GMG) algorithm is used, where 
for the multigrid level sub-solvers, a field-split (FS) preconditioner is proposed. 
The block structure of the FS preconditioner is derived using the physical variables as splitting strategy. 
To solve the subsystems originated by the FS preconditioning, an additive Schwarz (AS) block strategy is employed. 
The proposed field-split preconditioner is tested on biomedical FSI applications. 
Both 2D and 3D simulations are carried out considering aneurysm and venous valve geometries.
The performance of the FS preconditioner is compared with that of a second preconditioner of pure domain decomposition 
type. 
\end{abstract}
 

\section{Introduction}
Biomedical FSI problems are mathematically challenging because of the inherent non-linearity given by a domain
that moves as a function of the unknowns and because of the complex geometries involved 
(\cite{basting2017extended, bazilevs2010computational, crosetto2011fluid}). 
Numerical simulations of biomedical FSI problems often involve systems with a large number of unknowns
that make their solution computationally expensive in terms of time.
Moreover, the numerical stability could be compromised if forces are not computed correctly at the fluid-solid interface, or if 
large structural displacement arise, as in heart or venous valve simulations (\cite{causin2005added, le2001fluid}).

Many works on the numerical solution of FSI problems focus on the study of
preconditioning techniques for Krylov subspace solvers of block-structured linear systems 
(\cite{ipsen2001note, benzi2002preconditioning, benzi2005numerical, elman2011fast}).
Such a block structure arises, for instance, when finite element discretizations are adopted.
Multigrid algorithms, of either geometric or algebraic type, have been shown
to be reliable preconditioners for Krylov subspace algorithms (\cite{knoll1999multigrid, elman2001multigrid}). 
Multigrid procedures are appealing due to their potential for optimal computational
complexity and their ability to smooth out error components of both higher
and lower frequencies. Geometric multigrid (GMG) schemes have the advantage that projection 
and restriction operators are naturally obtained by the multilevel mesh discretization. On the other hand, 
algebraic multigrid (AMG) algorithms do not require the information of the multilevel mesh discretization but
they do not perform as well as perform than their geometric counterpart (\cite{watanabe2005comparison,volkov2016geometric}).

In this work, we consider a GMRES solver preconditioned with a GMG method,  
and propose a novel preconditioning technique for the multigrid level sub-solvers. 
For biomedical FSI simulations, fast algorithms are of fundamental importance, and an appropriate choice of 
the preconditioner guarantees that the system of equations is solved significantly faster. Hence,
the choice of preconditioners for the level sub-solvers can make a substantial difference 
on the overall performance of the solver.
The proposed preconditioner has a block structure that derives from using the physical variables as a splitting strategy,
therefore we refer to the proposed preconditioner as the field-split (FS) preconditioner. 
The FS preconditioner is a combination of 
physics-based (\cite{elman2008taxonomy}) and domain decomposition preconditioners (\cite{aulisa2015boundary}). 
These two strategies differ in the way the blocks are extracted from the original system matrix.  
Preconditioners constructed using a block structure associated with the
partial differential equation (PDE) system are categorized as physics-based.  
An example is the Schur complement preconditioner (\cite{bank1989class, murphy2000note, elman2002performance}).
Preconditioners based on a block structure associated with
the geometrical splitting of the domain are called of domain decomposition type. 
An example is the Vanka preconditioner (\cite{vanka1986block, manservisi2006numerical, wobker2009numerical}).
With the FS preconditioner, after the fields are separated with the proposed approach, 
the resulting subsystems are solved using an additive Schwarz (AS) block
strategy (\cite{cai1999restricted, brenner2008additive}). This is the reason why the FS preconditioner is a combination of the  
physics-based and domain decomposition approaches.
Preconditioners of field-split type have been used by the authors on Rayleigh-B\'{e}nard convection problems (\cite{ke2018new, Ke2016block, ke2018augmented}), and the intention of this paper is to investigate their capabilities for biomedical FSI problems.

The GMG is chosen over the AMG
because a monolithic Arbitrary Lagrangian-Eulerian (ALE) formulation is used for the FSI problem, where we exploit the 
high accuracy of the GMG inter-grid operators in preserving at each level the same fluid-solid decomposition
as at the finest level.

ALE schemes have been widely used to simulate biomedical FSI problems (\cite{turek2011numerical, molony2009fluid, bukavc2013fluid}).
A standard ALE approach moves the mesh to follow the elastic body movements.
Standard ALE algorithms with remeshing (\cite{margolin2003second, loubere2010reale}) and 
extended ALE schemes (\cite{basting2017extended, calandrini2018valve}) have been used to handle large structural displacements.
Other techniques proposed to solve 
FSI problems include fictitious domain methods (\cite{glowinski1994fictitious, baaijens2001fictitious, van2005three}), 
Lattice-Boltzmann methods (\cite{krafczyk1998analysis, fang2002lattice}), 
and the immersed boundary method (\cite{peskin2002immersed, peskin1989three}). 

For the coupling between the solid and fluid unknowns a monolithic approach has been chosen.
The coupling between the solid and fluid unknowns is a key factor to control stability issues. 
Monolithic schemes are, in general, more stable then partitioned schemes by means of solving simultaneously
for the fluid and solid part. In this way, techniques such as implicit discretization methods  
can be used, as well as strong coupled solvers for the whole system. 
Nevertheless, partitioned schemes have been investigated extensively in the literature (\cite{hou2012numerical})  
due to their relatively easier implementation.

\new{
Recent works in which GMRES has been used for the solution of a monolithic set of nonlinear equations are 
\cite{heil2004efficient, gee2011truly} and \cite{jodlbauer2019parallel}. 
In \cite{heil2004efficient}, preconditioners are developed based on block-triangular approximations of the Jacobian matrix, 
obtained by neglecting selected fluid-structure interaction blocks. 
In \cite{gee2011truly}, two  preconditioners based on algebraic multigrid techniques are applied to the Newton-Krylov solver. 
The first is based on a block Gauss-Seidel approach, where approximate inverses of the individual field blocks are based on an algebraic 
multigrid hierarchy. The second is based on a monolithic coarsening scheme for the coupled system. 
In \cite{jodlbauer2019parallel}, they construct preconditioners for Krylov subspace solvers based on block LDU-factorizations of the 
linearized FSI matrix.  
In \cite{richter2015monolithic}, GMG is used as a preconditioner for GMRES as we do, but for the smoothing strategy a partitioned iteration 
is used based on the idea provided in \cite{van2008space}. 
For an overview of some of the most popular methodologies to solve numerically the haemodynamic FSI systems, please refer to the introduction 
in \cite{crosetto2011parallel}, where a class of block triangular preconditioners is described obtained by exploiting the block-structure of 
the FSI linear system. 
}

To analyze the FS preconditioner for the level sub-solvers, we compare 
its performance with the preconditioner in \cite{FSIarticle}, which is of pure domain decomposition type.
In \cite{COUPLED2017,COMPDYN2017,calandrini2018magnetic} the Krylov solver in
\cite{FSIarticle} has been applied to solve FSI problems 
with biomedical applications, namely it has been used for stented cerebral aneurysm simulations in
\cite{COUPLED2017,COMPDYN2017} and for magnetic drug targeting (MDT) procedures in \cite{calandrini2018magnetic}. 
We intend to show that, 
a field-split approach is more suitable for biomedical FSI problems where a Krylov solver  
preconditioned with a multigrid method is adopted.

The biomedical applications considered for the tests include aneurysm and venous valve geometries. Both 2D and 3D simulations are carried out.
To further improve numerical stability, a Streamline Upwind Petrov Galerkin (SUPG) stabilization 
(\cite{franca1992stabilized, franca1993element}) 
is added to the momentum equation in our monolithic FSI formulation.

The paper is structured as follows. In Section \ref{modeling}, the weak monolithic formulation of the FSI
problem is described, and details are given on the SUPG stabilization technique used. In Section \ref{precond}, 
the field-split preconditioner for the level sub-solvers is described, together with the AS preconditioner
already presented in \cite{FSIarticle}. Moreover, the monolithic GMRES solver preconditioned with GMG is
briefly illustrated. In Section \ref{tests}, the numerical tests are presented.
Finally, in Section \ref{conclusions}, we draw our conclusions.

\section{Numerical modeling}\label{modeling}
This section describes the FSI formulation adopted in this paper. We present a monolithic ALE formulation assuming 
the fluid flow and the solid material to be incompressible. 
Since biomedical simulations are proposed in Section \ref{tests}, 
the solid structure refers to blood vessel walls, while the fluid considered is blood. 
To further improve stability, a Streamline Upwind Petrov Galerkin (SUPG) stabilization is employed. Following \cite{franca1993element}, the stabilization parameter, $\tau$, is computed by solving a generalized eigenvalue problem.
Before discussing the formulation of the problem, we introduce some notations.

\subsection{Notation}
For every time $ t \in [0,T]$,
let $\Omega^f_t \subset \mathbb{R}^n $ be an open set occupied only by a fluid,
and let $\Omega^s_t \subset \mathbb{R}^n $ be an open set occupied only by a solid.
In the following, any other symbol endowed with the superscripts $ f $ or $ s $ 
will refer to either the fluid or the solid part, respectively.
The boundary of the fluid and solid parts are denoted
as $ \partial \Omega^f_t $ and  $ \partial \Omega^s_t $, respectively. 
We define the parts of the boundary adjacent only to the fluid or only to the solid as $\Gamma^f_t$ and $\Gamma^s_t$,
such that $ \partial \Omega^f_t = \Gamma^f_t \cup \Gamma^i_t $ 
and $ \partial \Omega^s_t = \Gamma^s_t \cup \Gamma^i_t $, where $ \Gamma^i_t $ is the interface between fluid and solid.
The symbols $ \mathbf{n}^f $ and $ \mathbf{n}^s $ denote the outward unit normal fields defined on $ \partial \Omega^f_t$ and $ \partial \Omega^s_t $.
The open set $\Omega_t := \Omega^f_t \cup \Omega^s_t \cup \Gamma^i_t $ indicates
the current configuration of the overall physical domain.
The fluid and solid are immiscible, 
\begin{equation}
{\Omega^f_t} \cap {\Omega^s_t} = \emptyset \,,  
\end{equation}
and they interact through the nonempty interface $ \Gamma^i_t = \partial {\Omega^f_t} \cap \partial {\Omega^s_t} $.
The \textit{hat} notation is used to define $ \widehat{\Omega}^f := \Omega^f_0 $ and $ \widehat{\Omega}^s := \Omega^s_0 $.
Normally, they are referred to as the undeformed or stress-free configurations,
although the initial stresses
need not be identically zero either in the solid or in the fluid part.
Moreover, let us define $ \widehat{\Omega} := \Omega_0$ and $ \widehat{\Gamma}^i:= \Gamma^i_0 $.

The motion of the solid is followed in a Lagrangian way, therefore the domain $\widehat{\Omega}^s$ is a 
Lagrangian domain and it is initially occupied by the solid we observe. 
The domain $ \widehat{\Omega}^f $ is the domain on which we initially observe the fluid motion in an Eulerian way,
and it is called an ALE domain.
As a consequence of the solid movement, the domain on which we observe the fluid motion 
changes in time as well. Therefore, a deformation for both the fluid and solid domains needs to be defined.
The domain $ \Omega^f_t $ is occupied only by fluid at each time $ t $. 
The moving fluid or solid domains $ {\Omega^f_t} $ and $ {\Omega^s_t} $ are called Eulerian domains.

For the sake of brevity, the notations $\widehat{\nabla}$ or $\nabla$ refer to
the nabla symbolic operator in the fixed or moving frames, respectively.
In order to describe the motion of the fluid and solid domains, 
we define a $ t $-parametrized family of invertible and sufficiently regular mappings $ \mathcal{X}_t $,
called deformation mappings, given by a perturbation of the identity, so that
\begin{align}
 \mathcal{X}_t & : \widehat{\Omega} \rightarrow \Omega_t \,,  \quad \mathcal{X}_t(\widehat{\bm{x}})  := \widehat{\bm{x}} + {\bm{d}}(\widehat{\bm{x}},t)  \,, \label{mappings}
\end{align}
where the field $ {\bm{d}}(\widehat{\bm{x}},t) $ is the displacement field. 
These mappings are implicitly the solutions of the PDE system, so the definition given in \eqref{mappings} 
is actually a hypothesis on the regularity of the solutions of the PDE system.
For every $\widehat{\bm{x}} \in \widehat{\Omega}$ 
and $t \in [0,T]$, we also define
\begin{align}
  {\bm{F}} ({\bm{d}}(\widehat{\bm{x}},t) ) & = \widehat{\nabla}\mathcal{X}_t(\widehat{\bm{x}}) = I + \widehat{\nabla} {\bm{d}}(\widehat{\bm{x}},t) \,, \\ 
  {J}      ({\bm{d}}(\widehat{\bm{x}},t) ) & = \det {\bm{F}}({\bm{d}}(\widehat{\bm{x}},t) ) \,,\\
  {\bm{B}} ({\bm{d}}(\widehat{\bm{x}},t) ) & =  {\bm{F}} ({\bm{d}}(\widehat{\bm{x}},t) ) {\bm{F}}^\intercal 
  ({\bm{d}}(\widehat{\bm{x}},t) ) .
\end{align}
The symbols  $ {\bm{F}} $ and  $ {\bm{B}} $ denote 
the deformation gradient tensor and the left Cauchy-Green deformation tensor, respectively.

\subsection{Weak Monolithic Formulation}\label{weak}
In this section, the weak monolithic formulation of the FSI problem is described. 
With $\bm{d}$, $\bm{u}$ and $p$ we indicate displacement, velocity and pressure, respectively. 
Let us define the sets
 \begin{align}
 & \bm{V}^{\bm{f}}  :=  \bm{H}^1(\Omega^f_t) \,, \quad
   \bm{V}^{\bm{s}}  :=  \bm{H}^1(\Omega^s_t)
   \,, \\ 
 & \bm{V}^{\bm{f}}_0              :=  \bm{H}^1_{0}(\Omega^f_t;         {\Gamma}^f_{t,D} ) \;, 
   \bm{V}^{\bm{s}}_0              :=  \bm{H}^1_{0}(\Omega^s_t;         {\Gamma}^s_{t,D} ) \;, 
   {\bm{V}}^{\bm{f}}_{0,\bm{d}^f} :=  \bm{H}^1_{0}({\Omega}^f_t; {\Gamma}^f_{t,D,\bm{d}^f} )
   \,, \\ 
 & \bm{V} := 
  \{ \bm{v} = (\bm{v}^f,\bm{v}^s) \in 
  \bm{V}^{\bm{f}} \times \bm{V}^{\bm{s}} \text{ s. t. } \bm{v}^f = \bm{v}^s \text{ on } \Gamma^i_t \}  \,, \\
 & \bm{V}_0  := 
  \{ \bm{v} = (\bm{v}^f,\bm{v}^s) \in  
  \bm{V}^{\bm{f}}_0 \times \bm{V}^{\bm{s}}_0 \text{ s. t. } \bm{v}^f = \bm{v}^s \text{ on } \Gamma^i_t \} \,, 
 \end{align}
where $ {\Gamma}^f_{t,D} $ and $ {\Gamma}^f_{t,D,\bm{d}^f} $ are the subsets of $ {\Gamma}^f_t $
on which Dirichlet boundary conditions on the velocity $\bm{u}^f $ 
and on the fluid domain displacement $\bm{d}^f $ are enforced, respectively. 
Similarly, $ {\Gamma}^s_{t,D} $ denotes a subset of $ {\Gamma}^s_t $
on which Dirichlet boundary conditions on the displacement $\bm{d}^s $ are enforced. 
In the following, we will denote with the same symbol $ ( \cdot, \cdot )  $ 
the standard inner products either on $ L^2(\mathcal{O}) $, $ L^2(\mathcal{O})^n $
or $ L^2(\mathcal{O})^{n \times n} $, for any set $ \mathcal{O} \in \mathbb{R}^n$.

The monolithic weak FSI problem consists in finding 
$ (\bm{d}, \bm{u}, p) $ in 
$ \bm{V} \times \bm{V} \times L^2(\Omega_t) $ 
solution of a system that can be split into three parts:\\
the \textit{weak kinematic equations}
  \begin{align} 
  \left( {\bm{u}} - \frac{\partial {\bm{d}}}{\partial t}  , {\bm{\phi}}^{ks} \right)_{\widehat{\Omega}^s} 
   & = 0
  \qquad
  \forall \; {\bm{\phi}}^{ks} \in \bm{H}^1(\widehat{\Omega}^{s}), 
  \label{kinsld} \\  
   \left( k(\aled{\bm{x}})  \left(\aled{\nabla} {\bm{d}} + \big(\aled{\nabla} {\bm{d}}\big)^\intercal\right),  \aled{\nabla} {\bm{\phi}}^{kf}  \right)_{\widehat{\Omega}^f}
  & = 0 \qquad
 \forall \; {\bm{\phi}}^{kf} \in \bm{H}_0^1(\widehat{\Omega}^{f}; \widehat{\Gamma}^i) \cap \widehat{\bm{V}}^{\bm{f}}_{0,\bm{d}^f} \,, \label{kinfld}
  \end{align}
the \textit{weak monolithic momentum balance}
\begin{align}
 & \left( {\rho}^s  \frac{\partial {\bm{u}}}{\partial t} , {\bm{\phi}}^{m}       \right)_{{\Omega}_t^s}
 + \left(  {\bm{\sigma}}^s({\bm{d}},{p^s}) , {\nabla} {\bm{\phi}}^{m}  \right)_{{\Omega}_t^s} 
 - \left( {\rho}^s  {\bm{f}}^s    ,  {\bm{\phi}}^{m}      \right)_{{\Omega}_t^s}   
 \label{momentum_supg} \\
& + \left( \rho^f \frac{\partial \bm{u}}{\partial t} , \bm{\phi}^{m} \right)_{\Omega^f_t } 
  +  \left( \rho^f [ (\bm{u} - \frac{\partial \bm{d}}{\partial t} ) \cdot \nabla ] \bm{u} , \bm{\phi}^{m} \right)_{\Omega^f_t}
  \nonumber \\
& +  \left( \bm{\sigma}^f (\bm{u},p^f)     , \nabla \bm{\phi}^{m} \right)_{\Omega^f_t}
  -  \left( \rho^f \bm{f}^f   , \bm{\phi}^{m} \right)_{\Omega^f_t} \nonumber \\
& + \sum_{k=1}^{N^f} \left( \rho^f \frac{\partial \bm{u}}{\partial t} , \bm{\phi}_{SUPG}^{m} \right)_{T_t^k} 
  + \sum_{k=1}^{N^f}  \left( \rho^f [ (\bm{u} - \frac{\partial \bm{d}}{\partial t} ) \cdot \nabla ] \bm{u} , \bm{\phi}_{SUPG}^{m} \right)_{T_t^k} \nonumber \\
& - \sum_{k=1}^{N^f}  \left( \nabla \cdot \bm{\sigma}^f (\bm{u},p^f) , \bm{\phi}_{SUPG}^{m} \right)_{T_t^k}
  - \sum_{k=1}^{N^f}  \left( \rho^f \bm{f}^f  , \bm{\phi}_{SUPG}^{m} \right)_{T_t^k}
  = 0 \qquad
 \forall \; \bm{\phi}^{m} \in \bm{V}_0,\nonumber
\end{align}
and the \textit{weak mass continuity}
\begin{align}
\left( J( {\bm{d}} ) - 1 , {\phi}^{ps} \right)_{\widehat{\Omega}^s} 
 & = 0 
 \qquad
\forall \; {\phi}^{ps} \in L^2(\widehat{\Omega}^s),
\label{cntsld} \\ 
 \left( \nabla \cdot  \bm{u}  , \phi^{pf} \right)_{{\Omega}^f_t} 
 & = 0
\qquad
\forall \; \phi^{pf} \in L^2({\Omega}^f_t) \,.
\label{cntfld}
\end{align}
The above set of equations is known as a non-conservative ALE formulation.
In the following, we focus on the analysis of the momentum balance, and we briefly describe the role of the 
function $k(\aled{\bm{x}})$ in equation \eqref{kinfld}. 
For a full description of the PDE system \eqref{kinsld}-\eqref{cntfld}, please refer to \cite{FSIarticle}.
Equation \eqref{momentum_supg} describes, in a monolithic form, the solid and fluid momenta, which are also referred to as 
the incompressible non-linear elasticity equation and the Navier-Stokes equation, respectively.
The symbols $\rho^f$ and $\rho^s$ denote the mass densities for the fluid and solid part, respectively, whereas 
${{\mathbf{f}}}^f$  and ${{\mathbf{f}}}^s$ indicate the body force densities. 
The last two lines of equation \eqref{momentum_supg} describe the Streamline Upwind Petrov Galerkin (SUPG) method 
used to stabilize our FSI simulations. 
The set $\{T_{t}^k\}_{k=1}^{N^f}$ is a partition of $\overline{\Omega}_t^f$ such that $T_{t}^k \cap T_{t}^j$  has measure zero for $k\ne j$ and 
$\cup_{k=1}^{N^f} T_{t}^k = \overline{\Omega}^f_t$. 
Once a finite element discretization is introduced, 
$\{T_{t}^k\}_{k=1}^{N^f}$ will represent the finite element triangulation used.
The function $\bm{\phi}_{SUPG}^m$ is defined as
\begin{equation}
\bm{\phi}_{SUPG}^m=\tau \left( \rho^f [ (\bm{u} - \frac{\partial \bm{d}}{\partial t} ) \cdot \nabla ] \bm{\phi}^m \right)\,,
\end{equation}
where the description of the parameter $\tau$ is given below.
Notice that consistency between the fluid and solid region is automatically enforced with this stabilization technique. 
By definition of $\bm{\phi}_{SUPG}$, we have that $\bm{\phi}_{SUPG}=0$ in the solid, since 
$\bm u-\frac{\partial \bm d}{\partial t}=0$ on $\Omega^s_t$.
For the stability parameter $\tau$, the definition proposed by Franca and Madureira in \cite{franca1993element} is used, namely
\begin{equation}\label{tausupg}
 \tau=\dfrac{\xi(\mbox{Re}_k(\mathbf{x}))}{\sqrt{\lambda_k}|\bm{u}(\mathbf{x})|_2}\;,
\end{equation}
\begin{align}
 & \xi(\mbox{Re}_k(\mathbf{x}))=
\begin{cases}
  \mbox{Re}_k(\mathbf{x})\,,\;  \text{ if }\;\;\;  0\leq \mbox{Re}_k(\mathbf{x})<1\\
  1\,,\;\;\;\;\;\;\;\;\;\;  \text{ if }\;\;\;  \mbox{Re}_k(\mathbf{x})\geq 1
\end{cases}\;,
\end{align}
\begin{equation}
 \mbox{Re}_k(\mathbf{x})=\dfrac{|\bm{u}(\mathbf{x})|_2}{4\sqrt{\lambda_k}\nu(\mathbf{x})}\;,
\end{equation}
\begin{equation}\label{lambda_k}
 \lambda_k= \max_{\bm{u} \in \bm{V}_0(T^k_t)} \dfrac{||\nabla \cdot (\nabla \bm{u} + (\nabla \bm{u})^\intercal)||^2_{0,T^k_t}}{||(\nabla \bm{u} + (\nabla \bm{u})^\intercal)||^2_{0,T^k_t}}\;,
\end{equation}
\begin{align}
|\bm{u}(\mathbf{x})|_2= 
  \sum_{i=1}^{n}\big(|u_i(\mathbf{x})|^2\big)^{1/2}\,, \label{speed}
\end{align}
where $\nu(\mathbf{x})$ indicates the kinematic viscosity (ratio of the fluid viscosity to the fluid density).
The norm in equation \eqref{lambda_k} indicates the $L^2$ norm. 
With this design, no explicit computations of inverse estimate constants, nor the computation of mesh parameters is required.
The parameter $\lambda_k$ is computed as the largest eigenvalue of the following generalized eigenvalue problem:
for $k=1, \dots, N^f$ find $\bm{w} \in \bm{V}_0(T_t^k)$ and $\lambda_k$ such that 
\begin{equation}
(\nabla \cdot (\nabla \bm{w} +(\nabla \bm{w})^\intercal), \nabla \cdot (\nabla \bm{\phi} +(\nabla \bm{\phi})^\intercal))_{T_t^k} =
\lambda_k (\nabla \bm{w}, \nabla \bm{\phi})_{T_t^k},\;\;\forall \bm{\phi} \in \bm{V}_0(T^k_t).
\end{equation}
This problem is solved for the largest eigenvalue by the power method.

In equation \eqref{momentum_supg}, the interface physical condition of normal stress continuity 
is also enforced, where the boundary integrals disappear due to the condition
\begin{equation}
   \bm{\sigma}^s (\bm{d},p^s) \mathbf{n}^s + \bm{\sigma}^f (\bm{u},p^f) \mathbf{n}^f = 0 \quad \text{ on } {\Gamma}^i_t\,. \label{interface}
\end{equation}
For the solid stress tensor  $ {\bm{\sigma}}^s $, we consider incompressible Mooney-Rivlin, 
whose Lagrangian description is given for every $\widehat{\mathbf{x}} \in \widehat{\Omega}^s$ and $t \in [0,T]$ by
\begin{align}
{\bm{\sigma}}^s ({\bm{d}} , p^s )
   & = - p^s  \mathbf{I} + 2 C_1 {\mathbf{B}}({\bm{d}} ) -  2 C_2 ( {\mathbf{B}}
   ({\bm{d}} ) )^{-1} \;, \label{Mooney-Rivelin-cauchy}
\end{align}
where the constants $C_1$ and $C_2$ depend on the mechanical properties of the material. 
\new{For the numerical results presented in section \ref{tests}, we considered $C_1 = \frac{G}{3}$ and $C_2= \frac{C1}{2}$, 
where $G$ is the shear modulus. To compute $G$, we use the relation $E = 2G(1 + \nu)$, where $E$ is the Young's modulus
of the vessels wall and $\nu$ is the Poisson's ratio.}
Mooney-Rivlin constitutive relations have been used to represent the response of blood vessels (\cite{leach2010computational}). 
In \cite{tezduyar2007modelling}, it is shown that larger deformations of vessel walls can be achieved using the Mooney-Rivlin material 
instead of linearly elastic material models. 
For the fluid stress tensor $ {\bm{\sigma}}^f $, an incompressible Newtonian fluid is considered. 
In the literature, the assumption of Newtonian flow is generally accepted for numerical studies in large-sized arteries. 
For the treatment of stenotic vessels (\cite{guerciotti2018computational}) and cerebral aneurysms (\cite{cebral2005efficient}), 
it has been shown in that the computational results are moderately influenced by a non-Newtonian model.
For these reasons, blood is assumed as a Newtonian fluid in this paper, and its stress tensor is given 
for every ${\mathbf{x}} \in \Omega^f_t$ and $t \in [0,T]$ by
\begin{align}  \label{newtonian}
{\bm{\sigma}}^f ( {\bm{u}} , p^f ) = - p^f {\mathbf{I}} + \mu ( \nabla {\bm{u}} + (\nabla {\bm{u}})^\intercal ) \,,
\end{align}
where $\mu$ is the fluid viscosity.

In equation \eqref{kinfld}, the choice of the function $k(\aled{\bm{x}})$ is important for stability properties. 
This function controls the mesh deformation, and, according to the application considered, several choices of $k$ can be made.
If the element size is an issue for stability, a robust choice (\cite{tezduyar1992computation}) is 
\begin{equation}
 k(\aled{\bm{x}}) = \frac{1}{\mbox{V}_{el(\aled{\bm{x}})}},
\end{equation}
where $\mbox{V}_{el}$ is the volume of the mesh element that contains the $\aled{\bm{x}}$ coordinate.
In this way, a mesh with small elements where large displacements are expected,  
and large elements elsewhere is built.
In Section \ref{tests}, we often make the above choice for the function $k(\aled{\bm{x}})$.
However, for 2D simulations involving a venous valve geometry, the proximity of an element to the leaflets' tip affects stability, so  
the function $k(\aled{\bm{x}})$ is chosen to be a distance function 
that assures an homogeneous deformation throughout the entire mesh, and not only for those elements close to the valve leaflets. 
A more detailed discussion of the function $k(\aled{\bm{x}})$ adopted for the valve case can be found in Section \ref{valve_subsection}.

The monolithic weak FSI system \eqref{kinsld}-\eqref{cntfld} can be discretized in space using any appropriate choice of finite element spaces. 
Here, the finite element families chosen are biquadratic Lagrangian for both velocity and displacement, 
and \new{piecewise} linear discontinuous for pressure \new{(Q$_2$-P$_1$)}. 
For the time discretization, a second order Crank-Nicolson scheme is used for the momentum balance, and 
a backward finite difference method is adopted for the time derivatives in the integrands.
After the time and space discretizations, the resulting FSI system is linearized by an exact Newton linearization.
For more details about the fully discretized FSI system, 
please refer to \cite{FSIarticle}.

\section{Preconditioners}\label{precond}
The Krylov subspace solver for the linearized FSI system is GMRES preconditioned by a geometric multigrid algorithm, whose level sub-solvers are further preconditioned. 
In this section, the field-split preconditioner for the level sub-solvers is described together with a second preconditioner already presented in \cite{FSIarticle}. 
The proposed field-split preconditioner is a combination of the physics-based and domain decomposition approaches, while the second 
preconditioner is of pure domain decomposition type. 
Before discussing the structure of these two preconditioners, let us introduce some notation.
\subsection{Notation} 
A finite element discretization of the equations in Section \ref{weak}
leads to a system of algebraic equations.
In the following,  we use the symbols 
$K$ for the kinematic equation \eqref{kinsld} in the solid and on the fluid-solid interface,  
$A$ for the kinematic ALE displacement equation \eqref{kinfld} in the fluid,
$S$ for the momentum equation \eqref{momentum_supg} in the solid and on the fluid-solid interface, 
$F$ for the momentum equation \eqref{momentum_supg}  fluid respectively,
$V$ for the continuity equation \eqref{cntsld} in the solid, and
$W$ for the continuity equation \eqref{cntfld} in the fluid.
\subsection{The Structure of the Preconditioners} 
The ordering of the equations and variables determines the block structure of the Jacobian matrix $J$, 
resulting in a system that is more suitable to be solved using tailored preconditioners.
In the FSI literature, the most common way to rewrite a monolithic formulation is given 
by the following equation ordering: 
kinematic equations \eqref{kinsld}-\eqref{kinfld}, 
momentum equations \eqref{momentum_supg}, 
and continuity equations \eqref{cntsld}-\eqref{cntfld}. 
Thus,
\begin{equation} \label{jac_orig}
J = 
\mleft[
\begin{array} {c c | c c | c c}
K_{\bm{d}^s} 	&	0   		&	K_{\bm{u}^s}	&	0  		&	0 	    &   0 	\\
A_{\bm{d}^s}	&	A_{\bm{d}^f}	&	0 		&	0  		&  	0 	    &	0 	\\
\hline 
S_{\bm{d}^s}	&	S_{\bm{d}^f}	&	S_{\bm{u}^s}	&	S_{\bm{u}^f}	& S_{p^s} &	S_{p^f} 	\\
F_{\bm{d}^s}  	&       F_{\bm{d}^f}   	&	F_{\bm{u}^s}   	&   	F_{\bm{u}^f}	&	0 	 &	F_{p^f} \\
\hline 
V_{\bm{d}^s}		&	0 			& 	0 		& 	0 	&   	0     &  0	\\
W_{\bm{d}^s} 		&	W_{\bm{d}^f}	&   	W_{\bm{u}^s}	& 	W_{\bm{u}^f} 	&  0    & 0    
\end{array}
\mright],
\end{equation}
where the variables are ordered as $\bm{d}^s, \bm{d}^f, \bm{u}^s, \bm{u}^f, p^s, p^f$.
Note that in definition \eqref{jac_orig} there is no explicit reference to the interface.
The kinematic test functions on the interface are associated with the weak solid kinematic equation \eqref{kinsld} 
and are included in the first row of $J$. Recall that the test functions
in the weak fluid kinematic equation \eqref{kinfld} are in $\bm{H}_0^1(\widehat{\Omega}^{f}; \widehat{\Gamma}^i) \cap \widehat{\bm{V}}^{\bm{f}}_{0,\bm{d}^f}$, and they vanish on
the solid-fluid interface. Thus, the second row of $J$ does not directly affect the value of the displacement on the
interface. Similarly, for the weak momentum balance, test functions are divided into test functions for the solid part (third row) 
and test functions for the fluid part (fourth row), and contributions coming from nodes on the interface are accounted for in the solid part.

Recent works (\cite{aulisa2015multigrid, FSIarticle}) showed that the above equation ordering 
is ill-conditioned for steady-state problems, and it becomes unstable for
time-dependent problems as the time step increases. Such a behavior can be easily explained 
by looking at the diagonal terms $K_{\bm{d}^s}$ and
$S_{\bm{u}^s}$ in the first and third rows of \eqref{jac_orig}.
These terms correspond to the time derivatives of $\bm{d}^s$ and $\bm{u}^s$ in equations \eqref{kinsld} and \eqref{momentum_supg}, respectively. 
For steady state problems, such terms are identically zero, 
and for time dependent problems they become quite small for large time steps, resulting in a Jacobian matrix with zero or small diagonal terms. 
In \cite{aulisa2015multigrid, FSIarticle} a stable row/column pivoting alternative has been proposed, 
\begin{equation} 
J_1 = 
\mleft[
\begin{array} {c c | c c | c c}
S_{\bm{d}^s}		&	S_{\bm{d}^f}	&	S_{\bm{u}^s}	&	S_{\bm{u}^f}	& S_{p^s} &	S_{p^f} 	\\
A_{\bm{d}^s}		&	A_{\bm{d}^f}		&	0 	&	0  		&  	0 	    &	0 			\\
\hline 
K_{\bm{d}^s} 		&	0   	&	K_{\bm{u}^s}	&	0  			&	0 	    &   0 			\\
F_{\bm{d}^s}  	&       F_{\bm{d}^f}   	&	F_{\bm{u}^s}   	&   	F_{\bm{u}^f}		&	0 	    &	F_{p^f} 	\\
\hline 
V_{\bm{d}^s}		&	0 				& 	0 				& 	0 				&   	0  	    &  	0			\\
W_{\bm{d}^s} 		&	W_{\bm{d}^f}		&   	W_{\bm{u}^s }		& 	W_{\bm{u}^f} 		&  	0 	    & 	0    
\end{array}
\mright]. 
\end{equation}
For $J_1$, the ordering of the variables is $\bm{d}^s, \bm{d}^f, \bm{u}^s, \bm{u}^f, p^s, p^f$, as for $J$, and
the equations have been ordered in the following way: solid-interface momentum equations, 
fluid kinematic equation, solid-interface kinematic equations, fluid momentum equation,
solid continuity equation, and fluid continuity equation. In other words, $J_1$ is obtained from $J$ by permuting the first and the third row of equation \eqref{jac_orig}.
With this permutation, the time derivatives are moved off diagonal and are replaced by the matrix $S_{\bm{d}^s}$ and by
the lumped mass matrix $K_{\bm{u}^s}$. Note that the main contribution of $S_{\bm{d}^s}$ comes from the divergence of the stress tensor in equation  \eqref{momentum_supg}, which corresponds to the discretization of an elliptic operator. Therefore, this permutation assures a stable Jacobian matrix 
since $S_{\bm{d}^s}$ and $K_{\bm{u}^s}$ are always non-zero and somehow dominant.

A second stable row/column pivoting configuration, proposed here for the first time, is 
\begin{align} \label{matrix_J2}
J_2 &=
\mleft[
\begin{array}{c c c | c c c} 
S_{\bm{d}^s}	   &   S_{\bm{d}^f}	   & S_{p^s}  &  S_{\bm{u}^s}  &  S_{\bm{u}^f} 	    &  S_{p^f}  \\
A_{\bm{d}^s}	   &   A_{\bm{d}^f}      & 0 			&  0  			   &   0  			    &  0  \\                           												     
V_{\bm{d}^s}	   &   0    			   & 0 			&  0 			   &   0  			    &  0  \\
\hline 
K_{\bm{d}^s}     &   0  			   & 0 			&  K_{\bm{u}^s}  &   0  			    &  0  \\
F_{\bm{d}^s} 	 &   F_{\bm{d}^f}  & 0	&  F_{\bm{u}^s}  &   F_{\bm{u}^f }	    &  F_{p^f}   \\
W_{\bm{d}^s}	 &   W_{\bm{d}^f}  & 0	&  W_{\bm{u}^s } &   W_{\bm{u}^f} 	    &   0                                         															     
\end{array}
\mright].
\end{align}
For $J_2$, the variables are reordered as $\bm{d}^s, \bm{d}^f, p^s, \bm{u}^s, \bm{u}^f, p^f$, while
the equations are reordered in the following way: solid-interface momentum equations, 
fluid kinematic equation, solid continuity equation, 
solid-interface kinematic equations, fluid momentum equation, and fluid continuity equation. 
Comparing $J_1$ and $J_2$, the matrix $J_2$ is obtained from $J_1$ by permuting rows 3, 4 and 5 in 4, 5 and 3, and columns 
3, 4 and 5 in 4, 5 and 3.
This choice preserves the diagonal terms of $J_1$, i.e. $S_{\bm{d}^s}$, $A_{\bm{d}^f}$, $F_{\bm{u}^f }$ and $K_{\bm{u}^s}$, 
and it is more suited for the field-split preconditioner introduced below.

The systems associated with $J_1$ and $J_2$ are solved using 
GMRES preconditioned by GMG. 
At each multigrid level $l$ with $l \ge 1$, we use a damped Richardson iterative solver, whereas 
at the coarse level ($l=0$) a direct solver is employed.
Depending on the the block structure of the Jacobian matrix, the Richardson solver 
is preconditioned with an appropriate method.
The proposed additive field-split (FS) preconditioner is used on $J_2$, whereas the additive Schwarz (AS) preconditioner  
used for comparison is applied to $J_1$. 
The AS method is further preconditioned with an LU decomposition, 
while the FS preconditioner is further preconditioned with either AS-LU or AS-JACOBI.
For the sake of simplicity, we refer to these two combinations as the AS and the FS preconditioner, respectively. 
See Figure \ref{ASM} and \ref{FS} for a schematic representation of the two techniques. 
For both combinations, any type of multigrid cycle can be applied, 
such as V-cycle, F-cycle, or W-cycle, with any number of pre- and post-smoothing steps.
A version of the AS preconditioner was already developed and discussed in
\cite{aulisa2015multigrid, FSIarticle}, while the field-split preconditioner is the novelty of this paper. 
We remark that the FS preconditioner described in this work is an additive preconditioner. 
A multiplicative field-split preconditioner is not investigated since separate tests (not reported in this manuscript), 
show that it is not as efficient as the additive variant.

\begin{figure}[!h]
\centering
\begin{tikzpicture}[
node distance=3cm,
     terminal/.style={
      rectangle,rounded corners=3mm,minimum size=1cm,
       very thick,draw=black!100!blue,
       top color=white,
       }
]
{\node (gmres) [terminal] {GMRES}; }
{\node (gmg) [terminal,position=0:{0.6cm} from gmres] {GMG}; }
{\node (asm)  [terminal,position=0:{0.6cm} from gmg] { AS }; }
{\node (lu1)  [terminal,position=30:{2.0cm} from asm] {LU}; }
{\node (lu2)  [terminal,position=-30:{2.0cm} from asm] {LU}; }
{\path (gmres) edge[->,very thick] (gmg); }
{\path (gmg) edge[->,very thick] (asm); }
{\path (asm)  edge[->,very thick]    node [midway, above, sloped] (TextNode2) {$\bm{d}_j^s, \bm{u}_j^s, p_j^s$} node [midway, below, sloped] (TextNode2) {\small $j=1,\dots,n^s$} (lu1); }
{\path (asm)  edge[->,very thick]    node [midway, above, sloped] (TextNode2) {$\bm{d}_j^f, \bm{u}_j^f, p_j^f$} node [midway, below, sloped] (TextNode2) {\small $j=1,\dots,n^f$} (lu2); }
\end{tikzpicture}
\vspace{-8pt} \caption{$J_1$ solver and preconditioner cascade (AS).}
\label{ASM} 
\end{figure}
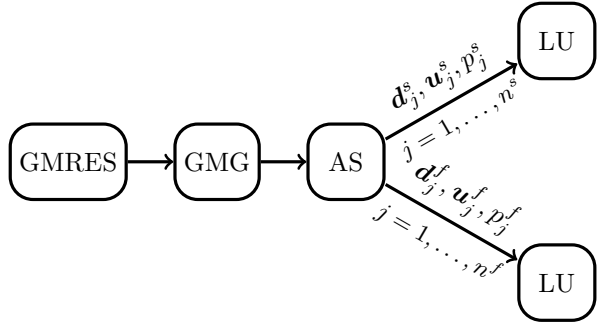

\begin{figure}[!h]
\centering
\begin{tikzpicture}[
node distance=3cm,
     terminal/.style={
      rectangle,rounded corners=3mm,minimum size=1cm,
       very thick,draw=black!100!blue,
       top color=white,
       }
]
{\node (gmres) [terminal] {GMRES}; }
{\node (gmg) [terminal,position=0:{0.6cm} from gmres] {GMG}; }
{\node (fs)  [terminal,position=0:{0.6cm} from gmg] { FS }; }
%
{\node (asm1) [terminal,position=35:{1.5cm} from fs] {AS$_1$}; }
{\node (lu1)   [terminal,position=20:{2.0cm} from asm1] { LU}; }
{\node (lu2)   [terminal,position=-20:{2.0cm} from asm1] { LU}; }
%
{\node (asm3)  [terminal,position=-35:{1.5cm} from fs] {AS$_2$ }; }
{\node (ilu)  [terminal,position=20:{1.0cm} from asm3] { JACOBI}; }
{\node (lu3)  [terminal,position=-20:{2.0cm} from asm3] { LU }; }
{\path (gmres) edge[->,very thick] (gmg); }
{\path (gmg) edge[->,very thick] (fs); }
{\path (fs) edge[->,very thick]   node [midway, above, sloped] (TextNode2) {$\bm{d},p^s$}  (asm1); }
{\path (asm1) edge[->,very thick] node [midway, above, sloped] (TextNode2) {$\bm{d}_j^s,p_j^s$} node [midway, below, sloped] (TextNode2) {\small $j=1,\dots,n^s$} (lu1); }
{\path (asm1) edge[->,very thick] node [midway, above, sloped] (TextNode2) {$\quad\bm{d}_j^f$} node [midway, below, sloped] (TextNode2) {\small $j=1,\dots,n^f$} (lu2); }
%
{\path (fs)  edge[->,very thick]   node [midway, above, sloped] (TextNode2) {$\bm{u},p^f$}  (asm3); }
{\path (asm3) edge[->,very thick]   node [midway, above, sloped] (TextNode2) {$\bm{u}^s$}  (ilu); }
{\path (asm3) edge[->,very thick]  node [midway, above, sloped] (TextNode2) {$\quad\bm{u}_j^f,p_j^f$} node [midway, below, sloped] (TextNode2) {\small $j=1,\dots,n^f$}  (lu3); }
\end{tikzpicture}
\vspace{-14pt} \caption{$J_2$ solver and preconditioner cascade (FS).}
\label{FS} 
\end{figure}
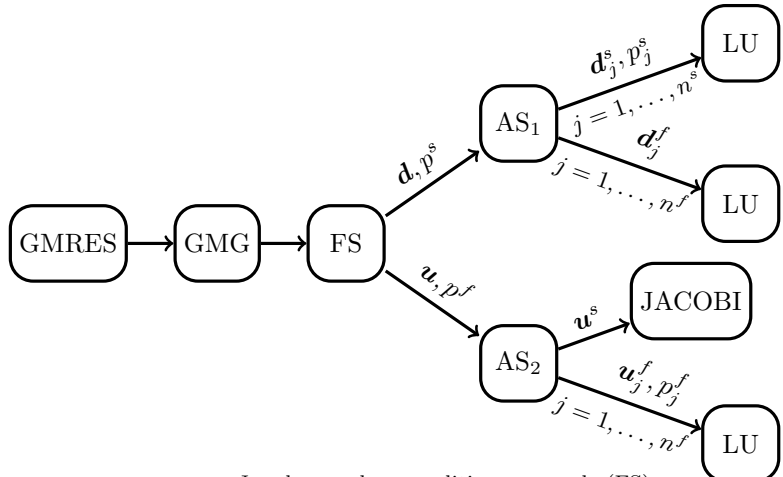

Now, let us discuss in detail the two preconditioning techniques, 
starting with the AS preconditioner, since an AS method is used in the novel FS preconditioner as well.
For the AS preconditioner, a locally multiplicative Vanka-type domain decomposition strategy is employed.
By locally multiplicative, we mean that the method is multiplicative within a single process, but additive among processes (\cite{15bcc4ae196041f1a0457f4a82ea9808}).
With the Vanka-type domain decomposition strategy, the physical domain is initially divided into a solid 
and a fluid sub-domain, and then each sub-domain is further divided into several
non-overlapping blocks, called Vanka blocks.
These non-overlapping Vanka blocks are such that their boundaries overlap with the boundaries of the mesh elements.
The variables within a single block are selected in the following way: on each block, all the pressure variables
whose support intersects the sub-domain are selected, and then all displacement and velocity variables
whose support intersects the support of the selected pressure are added.
Vanka-type domain decomposition is usually used for solving saddle point problems (\cite{vanka1986block}).
More details about this strategy can be found in \cite{aulisa2006computational, aulisa2015boundary, ke2018new}.
Depending on whether the Vanka block $j$ under consideration is in the fluid or in the solid subdomain, two different preconditioner matrices are created. The two matrices, $P^s_{1j}$ and $P^f_{1j}$, have the following form:
\begin{equation}
\mbox{AS}:\;
\left\{
\begin{array}{rcll} 
P^s_{1j} & = & 
\begin{bmatrix} 
S_{\bm{d}^s}	&	S_{\bm{u}^s}		& S_{p^s} \\
K_{\bm{d}^s} 	&	K_{\bm{u}^s}		&	0 	    \\
V_{\bm{d}^s}	& 	0 				&   	0  	    
\end{bmatrix}_{j_{_{_{_{}}}}}, & \mbox{ for } j=1,\dots,n^s
\\
P^f_{1j} & = & 
\begin{bmatrix} 
A_{\bm{d}^f}	&	0  		 &	0 			\\
F_{\bm{d}^f}    &   	F_{\bm{u}^f}	 &	F_{p^f} 	\\
W_{\bm{d}^f}	& 	W_{\bm{u}^f}	 & 	0    
\end{bmatrix}_j,
& \mbox{ for } j=1,\dots,n^f
\end{array} \right. .
\end{equation}
For the solid Vanka blocks, the preconditioner $P^s_{1j}$ is applied, and the kinematic, momentum and continuity equations are solved for $[ \bm{d}^s, \,\bm{u}^s,\,p^s]_j$ within each solid block. Similarly,  for the fluid Vanka blocks, the preconditioner $P^f_{1j}$ is applied, and the  kinematic, momentum and continuity equations are solved for $[\bm{d}^f,\, \bm{u}^f,\, p^f]_j$ within each fluid block. Note that the matrices forming $P^s_{1j}$ and $P^f_{1j}$ are relative small saddle-point matrices; therefore, the associated  saddle-point problems can be efficiently solved using LU decomposition. \new{A thorough explanation of how the AS block preconditioners $P^s_{1j}$ and $P^f_{1j}$ are applied to the Jacobian matrix $J_1$ can be found in \cite{smith2004domain}.}

\new{
In general, when seeking iteratively the solution of the system $J \bm{x} = \bm{f}$ for some given initial solution $\bm{x}^0$, in the 
AS case the solution at the $k^{th}+1$ iteration is obtained by
\begin{align}
&\quad \bm{y}^{0} = \bm{x}^k \,,\\
& \left.\begin{array}{l}
\bm{x}^{k+1}_j = \bm{x}^k_j + \left({P_{j}}\right)^{-1}(\bm{f} - J \bm{x}^k)_{j}  \\ \\
\bm{y}^{j}_i = \left\{
\begin{array}{l}
\bm{x}^{k+1}_{j_i} \mbox{ if the $i$ entry belongs to the block } j\\ \\
\bm{y}^{j-1}_i \mbox{ otherwise }
\end{array}
\right.
\end{array} \right\} \mbox{ for } j={1,\dots,n} \,, \\ 
&\quad \bm{x}^{k+1} = \bm{y}^{n}\,,
\end{align}
where the subspcript $j$ refers to the restriction to the $j^{th}$ subblock and the subscript $i$ refers to a single degree of freedom of the system.
In the additive case the block solution is updated using the solution available at the previous $k^{th}$ iteration.
Differently, in the multiplicative case the block solution is updated using the most recently updated block solution, namely
\begin{align}
&\bm{x}^{k+1}_j =  \bm{y}^{j-1} + \left({P_{j}}\right)^{-1}(\bm{f} - J \bm{y}^{j-1})_{j}, &\mbox{for } j={1,\dots,n}\,. 
\end{align}
The AS preconditioner available in PETSc (and used in this work) is multiplicative within a single process, but additive among processes.
The more processes are used the more the additive component dominates with respect to the multiplicative one, 
and in general additive is expected to be less performing than multiplicative.
}

The additive field-split preconditioner uses a Vanka-type domain decomposition strategy as well, but before that, the system 
is split according to the physical fields, i.e. $[\bm{d},\, p^s]$ and $[\bm{u},\, p^f]$.
This is the key feature of this novel preconditioner.
Based on the block structure of the matrix $J_2$, the Richardson iterative solver is preconditioned with the following matrix,
\begin{equation}
\mbox{FS}:\quad P_2 =
\mleft[
\begin{array} {ccc|ccc}
S_{\bm{d}^s}	   &   S_{\bm{d}^{f}}	 & S_{p^s}  &  0 			   &   0  			    &  0  		\\
A_{\bm{d}^s}	   &   A_{\bm{d}^f}      & 0       &  0  			   &   0  			    &  0  		\\            
V_{\bm{d}^s}	   &   0    			   & 0 			&  0 			   &   0  			    &  0  		\\
\hline 
0    			   &   0  			   & 0 			&  K_{\bm{u}^s}  &   0  			    &  0  \\
0 			   &   0 			   & 0 			&  F_{\bm{u}^s}  &   F_{\bm{u}^f }	&  F_{p^f}   \\
0			   &   0     			   & 0 			&  W_{\bm{u}^s }  &   W_{\bm{u}^f} 	&   0                             
\end{array}
\mright].
\end{equation}
The systems associated with the diagonal blocks are then preconditioned and solved separately using locally multiplicative AS block strategies 
denoted as AS$_1$ and AS$_2$, respectively,
\begin{equation}
\mbox{AS}_1: \;\left\{
\begin{array} {rcll}
{P_{2,1}^s}_j &= &
\begin{bmatrix} 
S_{\bm{d}^s}		&	S_{p^s}\\
V_{\bm{d}^s}		& 	0 			  	    
\end{bmatrix}_j, & \mbox{ for } j=1,\dots,n^s
\\
{P_{2,1}^f}_j &= &
\begin{bmatrix}
A_{\bm{d}^f}		
\end{bmatrix}_j,
& \mbox{ for } j=1,\dots,n^f
\end{array}\right.,
\end{equation}

\begin{equation}
\mbox{AS}_2:\;\left\{
\begin{array}{rcll}
{P_{2,2}^s} & = & 
\begin{bmatrix}
K_{\bm{u}^s}		
\end{bmatrix} &
\\
{P_{2,2}^f}_j & = & 
\begin{bmatrix}
   	F_{\bm{u}^f}		 &	F_{p^f} 	\\
 	W_{\bm{u}^f} 		 & 	0    
\end{bmatrix}_j,
& \mbox{ for } j=1,\dots,n^f
\end{array} \right. .
\end{equation}
The ${P_{2,1}^s}_j$ and ${P_{2,2}^f}_j$ blocks correspond to saddle point problems and are obtained using  
the Vanka-type domain decomposition strategy.
The ${P_{2,1}^f}_j$ blocks are obtained using an overlapping domain decomposition strategy, 
while ${P_{2,2}^s}$ corresponds to the lumped mass matrix and is not further decomposed.
The ${P_{2,1}^s}_j$, ${P_{2,1}^f}_j$ and ${P_{2,2}^f}_j$ blocks are relatively small matrices, 
so the systems associated with them can again be efficiently solved using LU decomposition. 
${P_{2,2}^s}$ is a diagonal matrix and the associate system can be solved exactly using one iteration of the Jacobi method.

The main difference between the AS and the FS preconditioner is that 
with the AS preconditioner the whole system of equations is solved together,
whereas with the FS preconditioner, the system is split into two smaller subsystems,  
and each is solved separately with its own AS strategy. 
The total number of sub-blocks doubles with FS where a domain decomposition strategy is applied, 
however the size of each sub-block is smaller than with the AS preconditioner. 
As the numerical examples in section \ref{tests} show, this size reduction results in less computational time, making
the FS preconditioner more efficient than the AS preconditioner, especially for three dimensional problems.

All of the above solvers and preconditioners have been implemented in FEMuS (\cite{femus-web-page}), 
an open-source finite element C++ library built on top of PETSc (\cite{petsc-web-page}).
Specifically, we have used the PCASM preconditioner from PETSc for the Vanka blocks 
and the PCFIELDSPLIT option for the splitting of the fields \cite{smith2004domain} . 
To perform the mesh partitioning among several processes, the METIS library (\cite{karypis1998software}) is used. 

\subsection{Geometric Multigrid Preconditioner}

A geometric multigrid preconditioner is used for the outer monolithic GMRES solver. 
In this section, we illustrate the multigrid prolongation and restriction operators 
associated with the Jacobian FSI matrices, $J_1$ and $J_2$.
For each level $l \in [0, N]$, let $h_l$ be the mesh size associated with the triangulation $\Omega_{h_l}$, where 
$\Omega_{h_l}$ is obtained recursively by midpoint refinement starting from a geometrically conforming coarse
triangulation $\Omega_{h_0}$.    
Let $\bm{\Phi}$ denote the piecewise biquadratic space and ${\Psi}$ the discontinuous piecewise linear space.    
Then, for each level $l$, ${\bm{\Phi}}({\Omega}_{h_{l}}) $ and $ {{\Psi}}({\Omega}_{h_{l}}) $ refer to
the finite element spaces associated with the triangulation $ \{\Omega_{h_l}\}_{l=1}^{N} $
with associated mesh size $ h_l $. 
With this choice of $\bm{\Phi}$ and ${\Psi}$, the level solution triplet $(\bm d, \bm u, p )_{h_l}$ belongs to the space  
$\bm{\Phi}({\Omega}_{h_{l}}) \times 
 \bm{\Phi}({\Omega}_{h_{l}}) \times 
 {\Psi}({\Omega}_{h_{l}}) $. 
We remark that, while continuity at the interface is imposed for displacement and velocity, the pressures 
$p^s$ and $p^f$ are unconstrained at the interface. For this reason the discontinuous piecewise linear space 
has been chosen for their discretization.

In the classical GMG for finite element methods, the prolongation operator 
\begin{align}
 \mathcal P_{l,l-1} & :
 \bm{\Phi}({\Omega}_{h_{l-1}}) \times 
 \bm{\Phi}({\Omega}_{h_{l-1}}) \times 
 {\Psi}({\Omega}_{h_{l-1}})
 \rightarrow 
 \bm{\Phi}({\Omega}_{h_{l}}) \times 
 \bm{\Phi}({\Omega}_{h_{l}}) \times 
 {\Psi}({\Omega}_{h_{l}})  
\end{align}
is the natural injection from a coarse to a fine space, 
and the restriction operator 
\begin{align}
  \mathcal R_{l-1,l} & :
 \bm{\Phi}({\Omega}_{h_{l}}) \times 
 \bm{\Phi}({\Omega}_{h_{l}}) \times 
 {\Psi}({\Omega}_{h_{l}})    
 \rightarrow 
 \bm{\Phi}({\Omega}_{h_{l-1}}) \times 
 \bm{\Phi}({\Omega}_{h_{l-1}}) \times 
 {\Psi}({\Omega}_{h_{l-1}})
\end{align}
is the adjoint of $ \mathcal P_{l,l-1}$ with respect to the $ L^2 $ inner product. Below, we show how
our restriction operator differs from the classical one. 

The matrix representations, $ P_{l,l-1}$ and $ R_{l-1,l} $, of the prolongation and restriction operators
depend on the row/column ordering of the Jacobian matrix, 
so the different structures of $J_1$ and $J_2$ generate two different pairs of matrices $ P_{l,l-1} $ and $ R_{l-1,l} $.
Namely, the matrix operators associated with $J_1$ have a block structure given by 
\begin{align*} 
& P_{1_{l,l-1}} =\left[
\begin{array}{cc | cc | cc}
   P_{\bm{d}^s,\bm{d}^s} &  0  &  0  &  0  &  0  & 0  \\
   P_{\bm{d}^f,\bm{d}^s} &  P_{\bm{d}^f,\bm{d}^f}  &  0  &  0  &  0  & 0  \\
\hline
   0  &  0  &  P_{\bm{u}^s,\bm{u}^s} & 0 & 0  & 0 \\
   0  &  0  &  P_{\bm{u}^f,\bm{u}^s} &  P_{\bm{u}^f, \bm{u}^f} & 0  & 0 \\
\hline
   0  &  0  &  0  &  0  &  P_{p^s,p^s} & 0    \\
   0  &  0  &  0  &  0  &     0        & P_{p^f,p^f} \\
\end{array} \right], \\
& R_{1_{l-1,l}} =\left[
\begin{array}{ cc | cc | cc}
  R_{S,S}    & 0 &  0  &  R_{S,F}      &   0&0\\
  0    &  R_{A,A}  &  0  &    0    &    0 &     0  \\
\hline
   0  &   0  &  R_{K,K }  &  0  &   0 & 0 \\
   0  &   0  &  0  &    R_{F ,F }  &   0          &     0      \\
\hline
   0             &         0            &    0   &     0              &     R_{V,V} & 0    \\
      0  & 0   &  0         &     0  
   &     0      &       R_{W,W} \\
\end{array} \right]
\,,
\end{align*}
whereas the matrix operators associated with $J_2$ have a block structure given by
\begin{align*} 
& P_{2_{l,l-1}} =\left[
\begin{array}{ccc|ccc}
   P_{\bm{d}^s,\bm{d}^s} &  0  &  0  &  0  &  0  & 0  \\
   P_{\bm{d}^f,\bm{d}^s} &  P_{\bm{d}^f,\bm{d}^f}  &  0  &  0  &  0  & 0  \\
    0  &  0  & P_{p^s,p^s} & 0  &  0  &   0    \\
\hline
   0 & 0  &  0  &  P_{\bm{u}^s,\bm{u}^s} & 0 & 0   \\
   0 & 0  &  0  &  P_{\bm{u}^f,\bm{u}^s} &  P_{\bm{u}^f, \bm{u}^f} & 0   \\
   0  &  0 & 0 &  0  &  0  &     P_{p^f,p^f}         \\
\end{array} \right], \\
& R_{2_{l-1,l}} =\left[
\begin{array}{ cc c |ccc}
  R_{S,S}  &0  & 0 &  0  &  R_{S,F}      &   0\\
  0    &  R_{A,A}  &  0  &    0    &    0 &     0  \\
   0             &         0            &    R_{V,V}   &     0              &     0 & 0    \\
\hline
 0&  0  &   0  &  R_{K,K }  &  0  &   0 \\
 0&  0  &   0  &  0  &    R_{F ,F }  &   0                \\
      0  & 0   &  0         &     0  
   &     0      &       R_{W,W} \\
\end{array} \right]
\,,
\end{align*}
where 
\begin{align*}
& R_{S,S} = R_{K,K}=(P_{\bm{d}^s,\bm{d}^s})^\intercal=(P_{\bm{u}^s,\bm{u}^s})^\intercal,
&& R_{A,A} = R_{F,F}=(P_{\bm{d}^f,\bm{d}^f})^\intercal=(P_{\bm{u}^f,\bm{u}^f})^\intercal,\\
& R_{V,V} =(P_{p^s,p^s})^\intercal, \quad R_{W,W} =(P_{p^f,p^f})^\intercal,
&& R_{S,F} =(P_{\bm{d}^f,\bm{d}^s})^\intercal. 
\end{align*}
Unlike the classical GMG, 
our restriction matrices are not simply the transposes of the associated prolongation matrices.
In particular, the blocks $P_{\bm{d}^f,\bm{d}^s}$ located at $P_{1_{l,l-1}}(2,1)$ and $P_{2_{l,l-1}}(2,1)$  
transpose into the blocks $R_{S,F}$ located at $R_{1_{l-1,l}}(1,4)$ and $R_{2_{l-1,l}}(1,5)$, respectively.
Such a structure for the matrix operators is due to the particular choice of row pivoting for the Jacobian matrices, $J_1$ and $J_2$.
Moreover, the transpose of the block $P_{\bm{u}^f,\bm{u}^s}$ in $P_{1_{l,l-1}}$ and $P_{2_{l,l-1}}$ 
does not appear in either restriction matrix, since 
the restriction of the kinematic equation in the solid-fluid interface 
does not receive any contribution from the fluid domain.

\section{Numerical results} \label{tests}
In this section, the performance of the proposed field-split preconditioner is 
tested on biomedical FSI problems. 
Aneurysm and venous valve geometries are considered, and 2D and 3D simulations are carried out. 
The performance of the proposed FS preconditioner 
is compared with that of the AS preconditioner, and 
\new{a comparison against a direct LU solver is also presented for the 2D geometries.  
To investigate the robustness of the schemes in terms of mesh refinement and time step size, a mesh independence and a time step
independence study are performed.}

All data are collected considering one single period, which we define as $1$ s. 
The time step considered is $\frac{1}{t_{step}}$ s, meaning that one period is composed of $t_{step}$ iterations.
For the aneurysm simulations $t_{step}=32$, while for the venous valve simulations $t_{step}=64$.
For every non-linear step $s$ at time step $t_i$, $i=1,\dots,t_{step}$, 
let us define the average convergence rate in the linear solvers 
as $\rho_{i,s} = \dfrac{r_{N_{i,s}}}{r_0}$, where $N_{i,s}$ is the number of linear steps 
(in the non-linear step $s$ at the time step $i$) 
and $r_0$ and $r_{N_{i,s}}$ represent the residuals of the initial and final linear iterations, respectively.
Over a single period, let us define the two quantities, $N$ and $\rho$, as 
\begin{align}
 N&=\dfrac{\sum_{i=1}^{t_{step}}\sum_{s=1}^{s_i^{max}} N_{i,s}}{\sum_{i=1}^{t_{step}}\sum_{s=1}^{s_i^{max}}1},\qquad
 \rho=\dfrac{\sum_{i=1}^{t_{step}}\sum_{s=1}^{s_i^{max}} \rho_{i,s}}{\sum_{i=1}^{t_{step}}\sum_{s=1}^{s_i^{max}}1},
\end{align}
where $s_i^{max}$ is the maximum number of non-linear steps at time step $t_i$, $i=1,\dots,t_{step}$.
The value $\rho$ indicates the average of the $\rho_{i,s}$'s over a single period, 
\new{i.e., the average convergence rate over a single period},
whereas $N$ represents the average of the $N_{i,s}$'s over a single period, 
\new{i.e., the average number of GMRES iterations over a single period}. 
With $s^{max}$ we indicate the average value of $s_i^{max}$ over a single period.

In the following tests, the values of $s^{max}$, $N$ and $\rho$ are computed for the FS 
and the AS preconditioner, \new{as well as for the LU  direct solver, when used. 
A comparison in terms of computational time is also presented.}  
The results for the 2D simulations are given first, followed by the 3D tests.


\subsection{2D Simulations}
Two 2D geometries are considered, involving a venous valve and a brain aneurysm. 
The AS preconditioner has already been tested on these two geometries in \cite{COUPLED2017, COMPDYN2017} and \cite{calandrini2018valve}.
A brief description of the function $k(\aled{\bm{x}})$ used in the venous valve simulations is also given. 

\subsubsection{Brain Aneurysm Geometry.}
\begin{figure}[ht]
\begin{center}
(a) \includegraphics[scale=0.45]{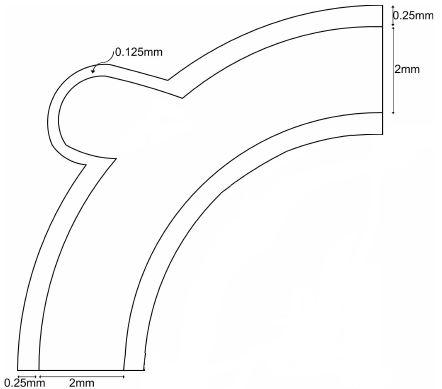}
(b) \includegraphics[scale=0.51]{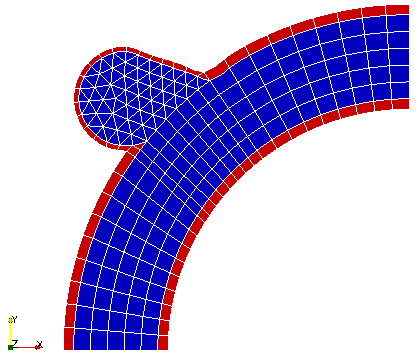} 
\caption{(a) 2D aneurysm geometry and lengths; (b) coarse mesh considered and material separation, where the red color indicates the solid domain and the blue color indicates the fluid domain.}\label{scheletro}
\end{center}
\end{figure}
The following 2D simulations of a cerebral aneurysm are based on a 2D hemodynamics model problem from \cite{turek2011numerical}. 
The geometry consists of a channel (lumen of the artery) of diameter $2$ mm with a wall thickness of $0.25$ mm (Figure \ref{scheletro} (a)). 
The aneurysm wall is typically thinner than that of
the healthy artery part; therefore the aneurysm has a wall thickness of $0.125$ mm.
For this geometry (and all the others as well) a hybrid mesh is employed. Quadrilateral elements are mainly used to mesh the lumen
and the arterial wall, while triangles are employed in the aneurysm bulge. The coarse mesh is displayed in Figure \ref{scheletro} (b), 
where the red color indicates the vessel wall and the blue color indicates the fluid domain. 
\new{The coarse mesh has $283$ elements and $1119$ biquadratic nodes, for a total of $5325$ degrees of freedom. 
In the simulations below, the coarse mesh is refined up to four times, and the last three levels are used for the quantitative analysis of the solver's performance.
These three levels have a total of $81492$,  $323524$ and $1289220$ degrees of freedom, respectively.}
The density and viscosity of blood are set to $1035$ kg/m$^3$ and $3.5 \times 10^{-3}$ Pa$\cdot$s, respectively.
For the elastic artery, a value of $1120$ kg/m$^3$ has been chosen for the density of the arterial wall, 
while the Young modulus and Poisson's ratio are set to $1.0$ MPa and $0.5$, respectively.
Inflow boundary conditions are specified as a pulsatile velocity profile moving from the right to the left part of the lumen
as in \cite{turek2011numerical} \begin{equation}
\bm{u}^f(0,y,t) = \left\langle (1 + 0.75\, \mbox{sin}(2\pi t)) u^f(0,y) , \, 0\right\rangle\;,
\end{equation}
where $u^f(0,y)$  is defined as parabolic inflow, namely \begin{equation}
u^f(0,y)=-0.05(1 - y^2),
\end{equation}
with $0 \le y\le R$, and $R=1$ mm is the lumen radius at the inlet.
A zero stress outflow boundary condition is considered at the lower left part of the artery.
The displacements are set equal to zero at the inlet and outlet of the artery, and
zero stress boundary conditions are imposed on the external aneurysm wall.

\subsubsection{Venous Valve Geometry.} \label{valve_subsection}
\begin{figure}
\begin{center}
 \includegraphics[scale=0.88]{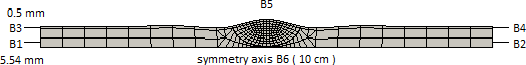}
 \caption{Vein and valve geometry with its coarse mesh, lengths and boundary names.}\label{valve_entire2}
\end{center}
\end{figure}
For the venous valve geometry, half of a blood vessel is considered,
namely only the motion of one of the two venous valve leaflets is analyzed. The other half of the
geometry and its mesh can be easily reconstructed by symmetry. 
Figure \ref{valve_entire2} shows the entire geometry, rotated of 90 degrees, together with its mesh,
lengths, and the labels given to the boundaries. In this geometry, we identify six different
boundaries: $B1$ and $B2$ represent the fluid boundaries at the bottom and top of the vein, respectively,
$B3$ and $B4$ represent the solid boundaries at the bottom and top of the vein, respectively, $B5$
indicates the lateral solid boundary, and $B6$ is the symmetry axis. The vein has a
lumen diameter of $5.54$ mm and a length of $10$ cm. The thickness of the vein wall is $0.5$ mm. For the
valve leaflet, a thickness of $0.065$ mm is specified. 
The elastic vein wall and the valve leaflet have the same density, $960$ kg/m$^3$, and Poisson's ratio, $0.5$.
The density and viscosity of blood are set to $1060$ kg/m$^3$ and $2.2 \times 10^{−3}$ Pa$\cdot$s, respectively.
At the fluid boundary $B1$, a normal stress boundary condition of the form
\begin{align}
\boldsymbol (\bm{\sigma}(\bm{u},p^f) \cdot \mathbf n) \cdot \mathbf n &= 15 \sin(2\pi t ) \quad \mbox{[Pa]}, \nonumber \\
\bm u  \cdot \boldsymbol \tau &= 0, \nonumber \\
\bm{d} &= \bm{0},
\end{align}
is specified, where $\bm{\tau}$ indicates the tangential vector to the boundary.
At the boundary $B2$, the same condition is applied but with the opposite sign. 
The vein is considered clamped; therefore at $B3$ and $B4$, the following condition is applied 
\begin{align}
 \bm{d} = \bm{0}.
\end{align}
At the lateral solid boundary $B5$, a stress boundary condition of the form
\begin{align}
 \bm{\sigma}(\bm{d},p^s) \cdot \mathbf n &= \bm{0},
\end{align}
is specified. Finally, at the symmetry axis $B6$, we require
\begin{equation}
\bm u  \cdot \mathbf n = 0,\qquad
\bm d  \cdot \mathbf n = 0, \qquad
\dfrac{\partial \mathbf u}{\partial \boldsymbol\tau} = 0, \qquad
\dfrac{\partial \mathbf d}{\partial \boldsymbol\tau} = 0. \nonumber
\end{equation}

\begin{figure}[t!]
 \begin{center}
  \includegraphics[scale=0.2, angle=-90]{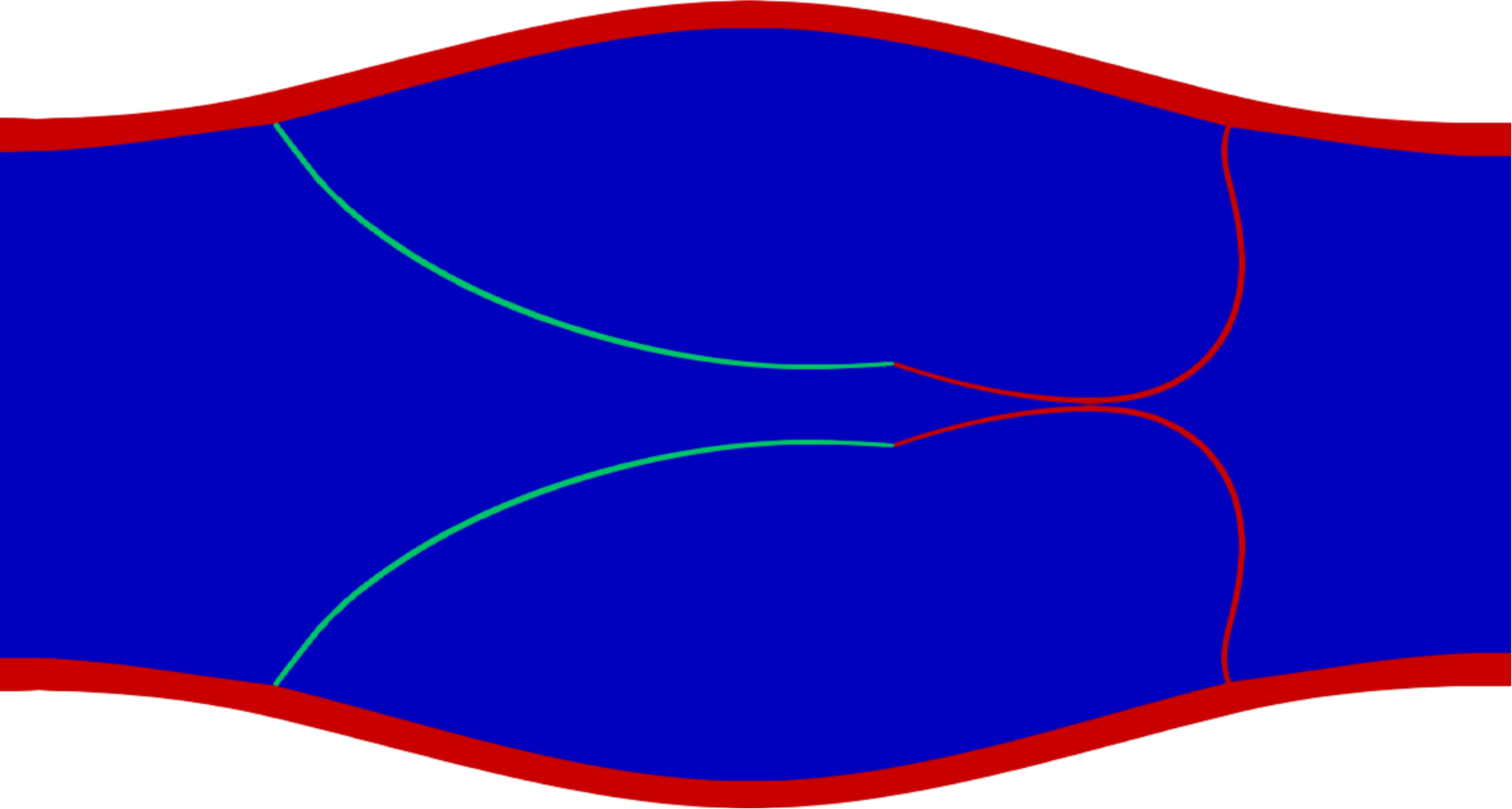}\hspace{1.1cm}
  \includegraphics[scale=0.2, angle=-90]{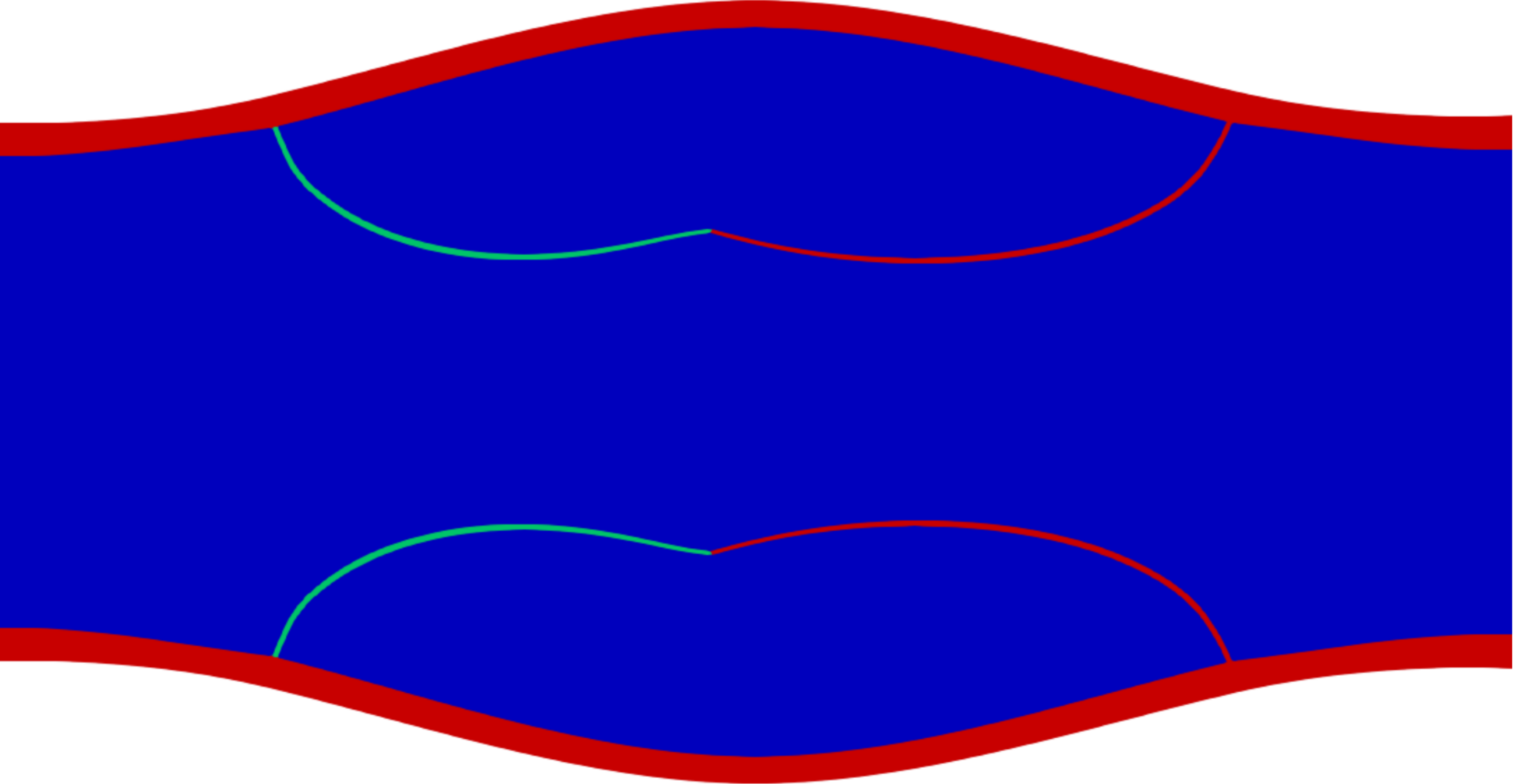} \\
  \end{center}
  \hspace{0.5cm}(a) \hspace{5.5cm}(b)
  \begin{center}
  \caption{Close (a) and open (b) leaflet configurations with highlighted \textit{fluid leaflet}, in green.}\label{zoom_n_tip}
 \end{center}
\end{figure}

\new{The valve coarse mesh, displayed in Figure \ref{valve_entire2}, has $292$ elements and $1245$ biquadratic nodes, for a total of $5856$ degrees of freedom.  
Similarly to the brain aneurysm case, the coarse mesh is refined up to four times and the last 3 levels are used for the quantitative analysis of the solver's performance.
These three levels have a total of $89604$,  $355716$ and $1417476$ degrees of freedom, respectively.}

For all the simulations involving this venous valve geometry, the function $k(\aled{\bm{x}})$ in equation \eqref{kinfld} 
is chosen as a distance function from a fluid element to the valve leaflet, to guarantee that there are no 
mesh entanglements in the fluid domain.
Given a fluid element in the mesh with center point $\aled{\bm{x}}$, the distance function $k(\aled{\bm{x}})$ is defined as
\begin{equation}
 k({\aled{\bm{x}}}) = \dfrac{a}{1+c\cdot k_1({\aled{\bm{x}}})}\,,
\end{equation}
where $a$, $c \in \mathbb{R}^{+}$, $a$, $c\ge 1$ and $k_1(\aled{\bm{x}})$ is a function that measures the distance from $\aled{\bm{x}}$, the center of the fluid element, to $m$, the middle point at the tip of the valve leaflet.
Thus, $k_1(\aled{\bm{{x}}})$ has the following expression
\begin{equation}
 k_1(\aled{\bm{x}}) = \sqrt{(\aled{x}_1 - m_1)^2+(\aled{x}_2 - m_2)^2},
\end{equation}
where the pair $(\aled{{x}}_1, \aled{{x}}_2)$  indicates the coordinates of $\aled{\bm{x}}$ and the pair $(m_1, m_2)$ indicates the coordinates of $m$. The value of $c$ considered in the numerical tests is $c = 10,000$, while for $a$ we consider either
$a=1$ or $a=100$. For the elements that belong to the \textit{fluid leaflet}, visible in Figures \ref{zoom_n_tip} (a) and (b) in green, we set $a=100$ to make these elements stiffer, whereas for the rest of the fluid elements, $a=1$. 
More details about the function $k(\aled{\bm{x}})$ can be found in \cite{calandrini2018valve}. 

\subsubsection{Results.}

\begin {table}[!tb]
\setlength\tabcolsep{5.5pt} 
\caption{2D Brain Aneurysm: mesh independence study for the solver}
\begin{center}
	\begin{tabular}{|c|c|c|c|c|c|c|c|c|} \hline	
	\multirow{2}{*} {}   & \multicolumn{8}{|c|}{2D Brain Aneurysm} \\ \cline {2-9}
	& \multicolumn{4}{|c|} {Additive Schwarz } & \multicolumn{4}{|c|} {Field-Split } \\ \hline		      
	 Level  & $s^{max}$ &  $N$  & $\rho$ & Time & $s^{max}$ & $N$   & $\rho$ & Time \\ \hline 
	     3  & 3.39      & 9.26  &  0.16  &  37.58s    &  3.39     & 9.12  &  0.15  &  37.84s   \\ \hline 
	     & \multicolumn{4}{|c|}{$\alpha = 0.94$} & \multicolumn{4}{|c|}{$\alpha=0.83$} \\ \hline
	     4  & 3.67      & 9.47  &  0.17  &  139.08s   &  3.67     & 8.86  &  0.15  &  119.39s    \\ \hline 
	     & \multicolumn{4}{|c|}{$\alpha = 0.97$} & \multicolumn{4}{|c|}{$\alpha=0.90$} \\ \hline
	     5  & 3.67      & 9.64  &  0.18  &  555.30s   &  3.67     & 9.11  &  0.16  &  454.00s    \\ \hline 
	\end{tabular}
\end{center}
\label{turek2D_mesh_indep}
\end{table}

\begin {table}[!tb]
\setlength\tabcolsep{5.5pt} 
\caption{2D Venous Valve: mesh independence study for the solver}
\begin{center}
	\begin{tabular}{|c|c|c|c|c|c|c|c|c|} \hline	
	\multirow{2}{*} {}   & \multicolumn{8}{|c|}{2D Venous Valve} \\ \cline {2-9}
	& \multicolumn{4}{|c|} {Additive Schwarz } & \multicolumn{4}{|c|} {Field-Split } \\ \hline		      
	 Level  & $s^{max}$ &  $N$  & $\rho$ & Time & $s^{max}$ & $N$   & $\rho$ & Time \\ \hline 
	     3  &   4.49    & 14.24 & 0.24   &  145.70s    &   4.49    & 14.33 & 0.24   &  139.72s   \\ \hline 
	     & \multicolumn{4}{|c|}{$\alpha = 1.00$} & \multicolumn{4}{|c|}{$\alpha=0.81$} \\ \hline
	     4  &   4.51    & 15.97 & 0.28   &  583.97s    &   4.51    & 16.21 & 0.28   &  429.58s   \\ \hline 
	     & \multicolumn{4}{|c|}{$\alpha = 1.03$} & \multicolumn{4}{|c|}{$\alpha=0.97$} \\ \hline
	     5  &   4.51    & 17.69 & 0.31   &  2508.16s   &   4.51    & 18.12 & 0.31   &  2039.69s   \\ \hline 
	\end{tabular}
\end{center}
\label{valve2D_mesh_indep}
\end{table}

\begin {table}[!tb]
\setlength\tabcolsep{5.5pt} 
\caption{2D Venous Valve: robustness of the solver in terms of time step size}
\begin{center}
	\begin{tabular}{|c|c|c|c|c|c|c|c|c|} \hline	
	\multirow{2}{*} {}   & \multicolumn{8}{|c|}{2D Venous Valve} \\ \cline {2-9}
	& \multicolumn{4}{|c|} {Additive Schwarz } & \multicolumn{4}{|c|} {Field-Split } \\ \hline		      
	 Time Step & $s^{max}$ &  $N$  & $\rho$ & Time & $s^{max}$ & $N$   & $\rho$ & Time \\ \hline 
	 1/64      & 4.51 & 15.97 & 0.28 & 583.97s  & 4.51 & 16.21 & 0.28 & 429.58s \\ \hline 
	 1/128     & 3.67 & 15.80 & 0.27 & 1024.95s & 3.67 & 16.11 & 0.28 & 767.03s \\ \hline
	 1/256     & 3.33 & 15.68 & 0.27 & 2739.24 & 3.33 & 16.14 & 0.28 & 1736.29s \\ \hline
	\end{tabular}
\end{center}
\label{valve2D_robustnessdt}
\end{table}

\new{First, a study to investigate the mesh independence of the proposed solvers is performed. 
Let the computational time be $$ T(\mbox{dofs}) = C \; \mbox{dofs}^\alpha, $$ for some constant $C$. Then, $\alpha=1$ is 
the optimal multigrid convergence rate for elliptic problems. If $\alpha < 1$ the convergence is better than linear and for $\alpha>1$ it is worse than linear.
To estimate the parameter $\alpha$ one can use the number of degrees of freedom and the computational times obtained using two different levels of refinement,
namely 
$$\alpha \approx \dfrac{\ln{ \dfrac{T(\mbox{dofs}_1)}{T(\mbox{dofs}_2)}} } {  \ln{ \dfrac{\mbox{dofs}_1} {\mbox{dofs}_2} } }.$$}

\new{Tables \ref{turek2D_mesh_indep} and 
\ref{valve2D_mesh_indep} show the values of $s^{max}$, $N$, $\rho$, $\alpha$ and the CPU time for the brain aneurysm and venous valve geometry, respectively. 
Three refinement levels are considered, namely 3, 4 and 5.
Sixteen processes were used to run these simulations in parallel. 
In Table \ref{turek2D_mesh_indep}, the values of $s^{max}$, $N$, $\rho$  are almost constant, whereas in Table \ref{valve2D_mesh_indep} 
$N$ and $\rho$ slightly increase with the number of uniform refinements. This increase, especially in the values of $\rho$,  
indicates that both preconditioning techniques are not completely independent from the mesh considered. A certain degree of dependence on the mesh 
is expected, since in hemodynamics simulations the meshes are constructed on complex geometries, and they may undergo severe deformation during the simulation, 
causing some deterioration in the convergence. This is especially visible in the second test, where in the closing phase the valve leaflet squeezes the fluid domain almost to nothing.
Nevertheless, the $\alpha$ parameter is always less than or equal to one. This is quite a remarkable result, given the nonlinearity and complexity of the FSI problem.
In terms of time, the FS preconditioner is faster than AS. For the most expensive simulations (5 refinement levels), the speed-up of FS over AS is 
18.24\% for the aneurysm case and 18.68\% for the venous valve. 
} 

\new{In terms of robustness with respect to the time step size, Table \ref{valve2D_robustnessdt} shows the results obtained for the 
more complex valve geometry. The values of $N$ and $\rho$ are stable for both solvers indicating that there is no significant change in the 
convergence rate. The value of $s^{max}$ decreases with the time step since, by reducing the time step size, the difference in the solution between two consecutive time intervals gets smaller.
Thus, a smaller number of non linear iterations is required for the solver to converge. 
In terms of computational time, FS performs better than AS, i.e. it is less expensive than AS for all time step sizes considered.} 

\new{In terms of scalability of the solvers with respect to the number of processes, Tables \ref{turek2D_ksp&time} and \ref{valve2D_ksp&time} 
show the results given for the brain aneurysm and venous valve geometry, respectively. 
The number of processes used is 4, 8 and 16, and for the two geometries the fourth uniform refinement level is considered. 
To have a more complete picture, a direct LU solver is added for comparison. This third solver 
is referred to as MLU, since we used the implementation available in the MUMPS library \cite{MUMPS}.} 
\new{From Table \ref{turek2D_ksp&time}, we can see that the values of $s^{max}$ and $\rho$ obtained with the AS and FS preconditioners are 
very similar. Slightly bigger differences can be found in the values of $N$. An analogous situation is present in Table \ref{valve2D_ksp&time}. 
This means that the two solvers exhibit similar convergence properties. Moreover, the convergence does not deteriorate increasing the number of 
processes. 
Looking at the results obtained with MLU, better convergence properties are observed,  
although the computational time is larger.} 

\new{The solver time decreases as the number of processes increases. 
For the aneurysm geometry, we have that when using 8 processes instead of 4, the reduction in the solver time is 43.96\% with the AS preconditioner, 
and 41.78\% with the FS preconditioner. With 16 processes instead of 8, the reduction is 40.48\% with the AS preconditioner, 
and 37.31\% with the FS preconditioner. 
For the 2D valve case, we have that, when using 8 processes instead of 4, the reduction in the solver time is 46.47\% with the AS preconditioner, 
and 47.03\% with the FS preconditioner; whereas using 16 processes instead of 8, the reduction is 37.58\%
with the AS preconditioner, and 40.47\% with the FS preconditioner.
Therefore, the performances of the two preconditioning techniques are comparable in terms of process scalability, and the scalability 
is not significantly impacted by the geometry. In an ideal situation, we would like the time reduction to be exactly half (50\%), when doubling the number of 
processes. Such a scenario is not achieved here, due to the use of the Schwartz domain decomposition preconditioner available in PETSc. We remind that this preconditioner is multiplicative within a process 
but additive among processes. The more processes are used the more the additive component dominates with respect to the multiplicative one, 
and in general additive is expected to be less performing than multiplicative.}

\new{The MLU solver has an even poorer performance in terms of 
process scalability, in fact, for the aneurysm case, the time reduction is only 22.46\% from 4 to 8 processes and 15.18\% from 8 to 16. Similarly, 
for the valve case, the time reduction is just 22.29\% from 4 to 8 processes and 18.60\% from 8 to 16.} 

\new{
Comparing the computational times, the FS preconditioner is considerably faster than the AS and MLU for all cases considered.
From Table \ref{turek2D_ksp&time} we see that using the FS preconditioner, the time is reduced by 21.55\% over AS and 36.80\% over MLU with 4 processes, 
by 18.49\% over AS and 62.78\% over MLU with 8 processes, and by 14.16\% over AS and 72.49\% over MLU with 16 processes. 
Hence, there is an average reduction of time of 18.07\% over AS and 57.37\% over MLU. 
Focusing on Table \ref{valve2D_ksp&time}, the time is reduced by 22.07\% over AS and 29.46\% over MLU with 4 processes, 
by 22.85\% over AS and 51.91\% over MLU with 8 processes, and by 26.44\% over AS and 64.83\% over MLU with 16 processes. 
Hence, for the valve case, there is an average reduction of 23.77\% over AS and 48.73\% over MLU.}

\begin{figure}[!tb]
  \begin{center}
     (a) \includegraphics[scale=0.54]{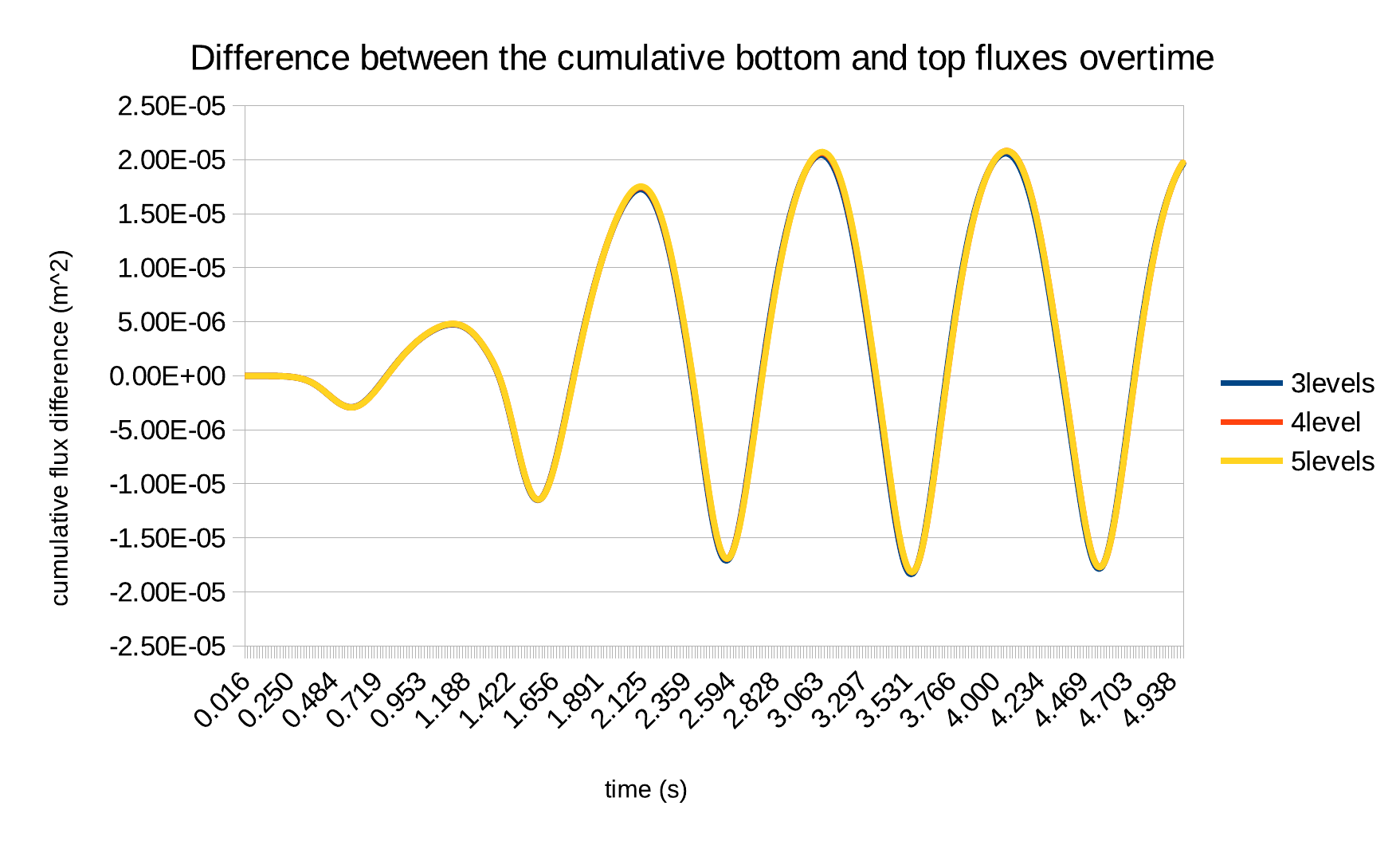}\\
     \hspace{0.7cm} (b) \hspace{-0.2cm} \includegraphics[scale=0.63]{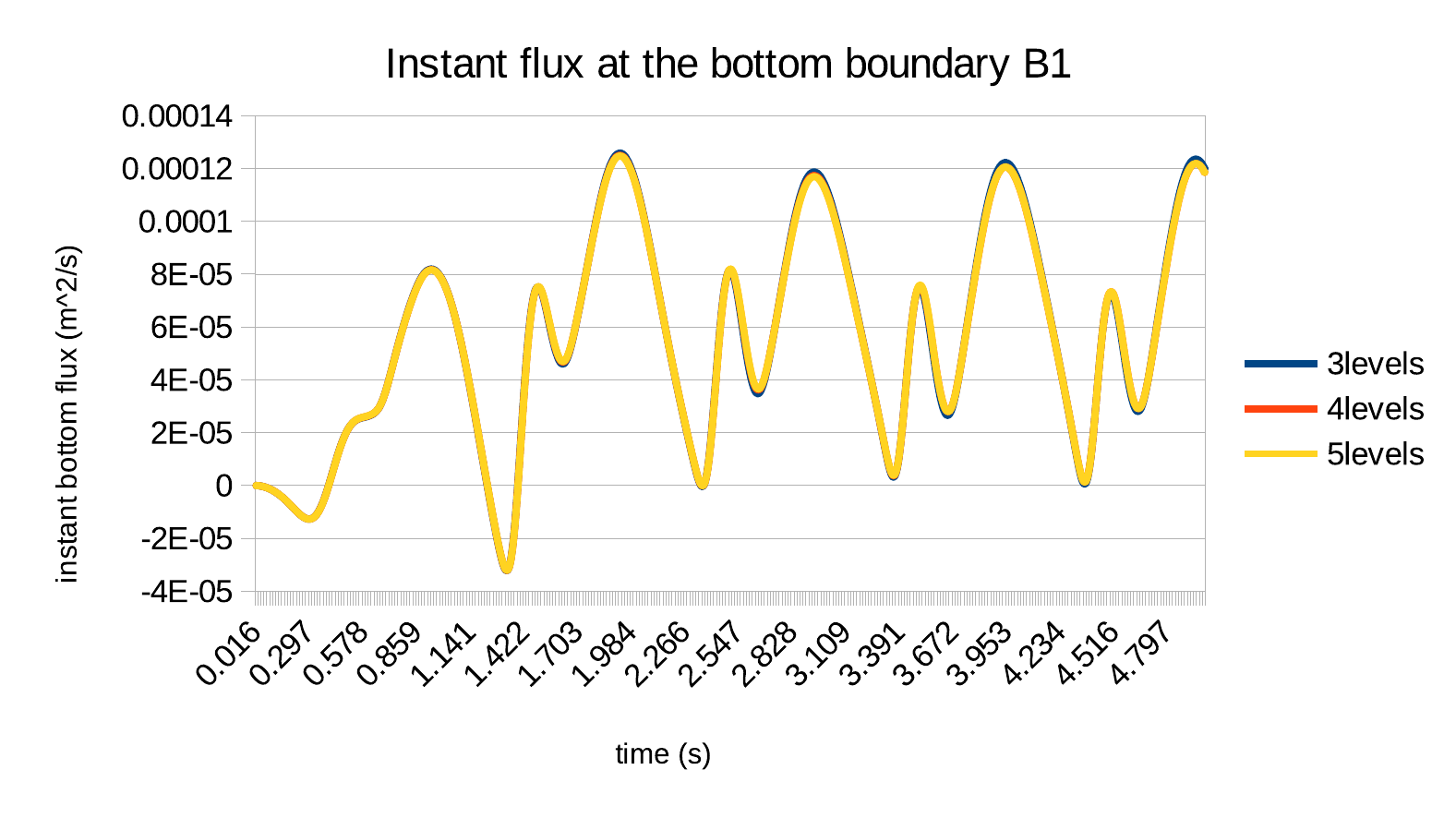}
    \caption{(a) Difference over time between the cumulative bottom and top fluxes, $Q_1(t) - Q_2(t)$; (b) instant flux $q_1(x)$ at the bottom boundary $B_1$.}\label{fluxes}
  \end{center}
\end{figure}

%
\new{Finally, we show that with the proposed FSI modeling, our solvers are able to consistently capture physical quantities of
interest. The results obtained with the AS or FS preconditioner are identical, so only those obtained using FS are shown. 
Focusing on the 2D valve geometry, Figure \ref{fluxes} shows the difference between the cumulative bottom ($B_1$) and top ($B_2$) 
fluxes over time and the instant flux at the bottom boundary $B_1$. 
These two quantities of interest are computed for three different levels of mesh refinements, namely levels 3, 4 and 5. 
Total flux (instant flux) at a boundary $B_i$, $i=1, 2$, is defined as \begin{equation}
q_i(t)= \int_{B_i} \bm{u}(x,t) \cdot \mathbf{n} \; ds,                                                             
\end{equation}
where $\mathbf{n}$ is the inlet normal in $B_1$ and the outer normal in $B_2$.
The cumulative total flux over time is the integral over time of the total flux, namely \begin{equation}
Q_i(t)= \int_{0}^t q_i(\tau) \; d\tau.                                                             
\end{equation}
As expected, the difference between the cumulative bottom and top fluxes over time always oscillates around 0 for all three simulations, 
meaning that the fluxes at the boundaries $B_1$ and $B_2$ of the vein correctly balance each other over time. The discrepancies in the solutions 
for the three levels of refinements are negligible, therefore we conclude that the solution is correctly represented.
The instant flux at the bottom of the vein geometry shows that some backflow occurs at the inlet. This natural phenomenon occurs because, 
when the valve closes, it undergoes a large deformation (especially for a small Young's modulus $E$), and blood may be pushed down, 
causing an exiting flux at the inlet. Again, the differences in the solutions for the three levels of refinements are negligible.}  

\begin {table}[!tb]
\setlength\tabcolsep{3.2pt} 
\caption{2D Brain Aneurysm: solver tests - 4 refinements}
\begin{center}
	\begin{tabular}{|c|c|c|c|c|c|c|c|c|c|c|c|c|} \hline	
	\multirow{2}{*} {}   & \multicolumn{12}{|c|}{2D Brain Aneurysm} \\ \cline {2-13}
	& \multicolumn{4}{|c|} {Additive Schwarz } & \multicolumn{4}{|c|} {Field-Split } & \multicolumn{4}{|c|} {MLU} \\ \hline		      
	 procs  & $s^{max}$   & $N$  & $\rho$ & Time     & $s^{max}$ & $N$  & $\rho$ & Time     & $s^{max}$ & $N$      & $\rho$    & Time \\ \hline 
	     4  & 3.67 & 9.59 &  0.18  &  417.02s &  3.67     & 8.90 &  0.16  & 327.13 s & 3.67      & 1 & 6.27E-011 & 659.89s \\ \hline 
	     8  & 3.67 & 10.15&  0.19  &  233.68s &  3.67     & 9.32 &  0.17  & 190.45 s & 3.67      & 1 & 6.86E-011 & 511.66s \\ \hline 
	    16  & 3.67 & 9.47 &  0.17  &  139.08s &  3.67     & 8.86 &  0.15  & 119.39 s & 3.67      & 1 & 7.02E-011 & 433.96s \\ \hline 
	\end{tabular}
\end{center}
\label{turek2D_ksp&time}
\end{table}

\begin {table}[!tb]
\setlength\tabcolsep{2.8pt} 
\caption{2D Venous Valve: solver tests - 4 refinements}
\begin{center}
	\begin{tabular}{|c|c|c|c|c|c|c|c|c|c|c|c|c|} \hline	
	\multirow{2}{*} {}   & \multicolumn{12}{|c|}{2D Venous Valve} \\ \cline {2-13}
	& \multicolumn{4}{|c|} {Additive Schwarz } & \multicolumn{4}{|c|} {Field-Split } & \multicolumn{4}{|c|} {MLU} \\ \hline		      
	 procs  & $s^{max}$ &  $N$  & $\rho$ & Time     & $s^{max}$ & $N$   & $\rho$ & Time     & $s^{max}$ & $N$      & $\rho$ & Time\\ \hline 
	     4  &   4.51    & 16.48 & 0.29   & 1748.17s &   4.51    & 16.60 & 0.29   & 1362.29s & 4.51      & 1 & 1.31E-08 & 1931.23s  \\ \hline 
	     8  &   4.51    & 15.84 & 0.28   & 935.50s  &   4.51    & 16.55 & 0.29   & 721.65s  & 4.51      & 1 & 1.01E-08 & 1500.80s \\ \hline 
	    16  &   4.51    & 15.97 & 0.28   & 583.97s  &   4.51    & 16.21 & 0.28   & 429.58s  & 4.51      & 1 & 9.73E-09 & 1221.58s \\ \hline 
	\end{tabular}
\end{center}
\label{valve2D_ksp&time}
\end{table}

\subsection{3D Simulations}
Two 3D geometries are presented, a brain aneurysm and an aortic aneurysm.
The AS preconditioner has already been tested on the 3D brain aneurysm geometry in \cite{COUPLED2017} and \cite{COMPDYN2017}.

\subsubsection{Brain Aneurysm Geometry.}
\begin{figure}[ht]
  \begin{center}
    (a) \includegraphics[scale=0.53]{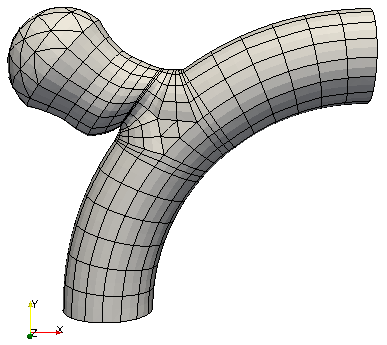}\hspace{0.1cm}
    (b) \includegraphics[scale=0.55]{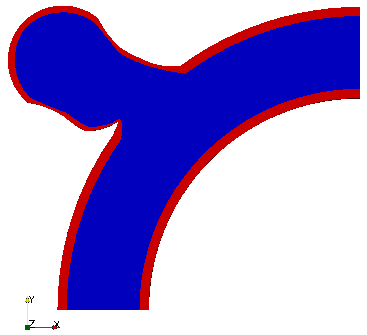}
    \caption{(a) 3D brain aneurysm geometry and mesh; (b) cross-section of the 3D geometry where the red color indicates the vessel wall and the blue color indicates the fluid domain.}\label{figure3D}
  \end{center}
\end{figure}
The cerebral aneurysm considered here is a 3D extension of the 2D geometry analyzed above.
Now, the aneurysm wall is assumed to have a uniform thickness equal to $0.25$ mm. Once again, a hybrid mesh is employed.
Wedges are utilized to mesh the artery lumen, whereas hexahedra are used for the arterial wall. Tetrahedral elements are mainly
employed in the aneurysm cavity. Figure \ref{figure3D} (a) shows the geometry with its mesh, whereas Figure 
\ref{figure3D} (b) is a cross-section of the 3D geometry where the red color indicates the vessel wall and the blue color indicates
the fluid domain. 
\new{The mesh has $268$ elements and $2169$ triquadratic nodes, thus the number of degrees of freedom is $14086$. 
In the simulations performed below, this mesh is refined up to two times, 
and in this case we have $17152$ elements, $121253$ triquadratic nodes and so $796126$ degrees of freedom.}
The same physical parameters used in the 2D case are employed for both the artery wall and blood.
The inflow boundary condition is a pulsatile velocity profile moving from the right to the left part of the lumen
described by 
\begin{equation}
\bm{u}^f(0,y,z,t) = \left\langle -0.3(1 - r^2)(1 + 0.75\, \mbox{sin}(2\pi t)),0,0\right\rangle\;,
\end{equation}
where $r=\sqrt{y^2+z^2}$, with $0\le r\le R$, and $R$ is the lumen radius at the inlet.
At the outlet, a pressure condition of the form $p = 0$ has been imposed.
The boundary displacements at the inlet and outlet of the artery are set to zero.

\subsubsection{Aortic Aneurysm Geometry.}
\begin{figure}[t!]
 \begin{center}
  (a) \includegraphics[scale=0.60]{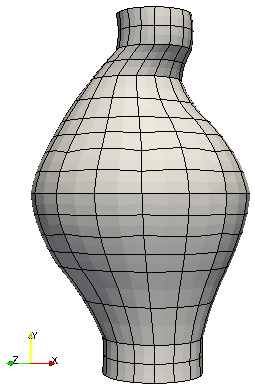} \hspace{1.5cm}
  (b) \includegraphics[scale=0.64]{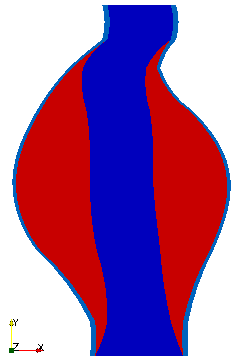}
  \caption{(a) aortic aneurysm geometry and mesh, (b) cross-section of the aortic aneurysm geometry, where the blue color shows the lumen, the red color shows the clotted blood and the light blue colors shows the artery wall.}\label{aortic}
 \end{center}
\end{figure}

The aortic aneurysm geometry and mesh are shown in Figure \ref{aortic} (a).
\new{The mesh has $208$ elements and $1715$ triquadratic nodes, thus the number of degrees of freedom is $11122$. 
In the simulations performed below, this mesh is refined up to two times, 
and in this case we have $13312$ elements, $98537$ triquadratic nodes and so $644470$ degrees of freedom.}
The arterial wall has a thickness of $1.5$ mm at the inlet (bottom part of Figure \ref{aortic} (a)) and outlet (top part of Figure \ref{aortic} (a)) and it gets thinner getting closer to the aneurysm, reaching a 1 mm thickness.
The lumen diameter at the inlet is $2.5$ cm, while at the outlet it reduces to $1.8$ cm.
The maximum aneurysm diameter is $6$ cm. The hybrid mesh employed uses hexahedral elements for the arterial wall and wedges and tetrahedral elements
for the lumen. For both the artery wall and blood, the same physical parameters used in the 2D and 3D brain aneurysm cases are employed. 
Blood flow often slows in the bulging section of an aortic aneurysm, causing clots to form.
Therefore, a second very soft solid is introduced in the aneurysm cavity to simulate clotted blood.
The Young's modulus of this second solid is $3\times10^{-3}$ MPa.
Figure \ref{aortic} (b) shows a cross-section of the geometry where the lumen, clotted blood and the artery wall are visible. 
The inflow boundary condition is a pulsatile velocity profile moving from the bottom to the top part of the lumen described by 
\begin{equation}
\bm{u}^f(x,0,z,t) = \left\langle 0,-\dfrac{0.01}{0.81}(0.81 - r^2)(1 + 0.75\, \mbox{sin}(2\pi t)),0\right\rangle\;,
\end{equation}
where $r=\sqrt{x^2+z^2}$, with $0\le r\le R$, and $R$ is the lumen radius at the inlet.
At the outlet, a time-dependent normal stress boundary condition of the form 
$$ \boldsymbol (\bm{\sigma}(\bm{u},p^f) \cdot \mathbf n) \cdot \mathbf n= 12500 + 2500 \sin(2\pi t), \quad \boldsymbol (\bm{\sigma}(\bm{u},p^f) 
\cdot \mathbf n) \cdot \boldsymbol\tau= 0$$ 
has been imposed, 
to simulate the physiological pressure range 80 - 120 mm Hg.
The artery at both the inlet and the outlet is considered clamped.

\subsubsection{Results.}

For both geometries, 3 uniform refinement levels have been considered. 
\new{
Tables \ref{turek3D_ksp&time} and \ref{AAA3D_ksp&time} show the results for the 3D brain aneurysm and the aortic aneurysm geometry, 
respectively. For both geometries, the values of $s^{max}$, $N$ and $\rho$ obtained with the two preconditioners
are comparable, so the two techniques have similar convergence properties. 
The values of $N$ and $\rho$ increase with the number of processes for both solvers, meaning that the convergence is negatively affected 
by increasing the number of processes.
This is somewhat different from the 2D results, where convergence was better preserved, but 3D problems are stiffer than 2D ones, so 
it is reasonable to expect a worsening in the converge rate and in the number of linear iterations.}
In terms of process scalability, we have that the performances of the two preconditioning techniques are comparable, as in the 2D cases. 
As expected, the solver time decreases as the number of processes increases.
For the brain aneurysm, using 8 processes instead of 4, the reduction in the solver time is 39.56\%
with the AS preconditioner and 33.60\% with the FS preconditioner; whereas using 16
processes instead of 8, the reduction is 29.67\% with the AS preconditioner and 29.66\%
with the FS preconditioner. 
For the aortic aneurysm, using 8 processes instead of 4, the reduction in the solver time is 20.41\% 
with the AS preconditioner and 14.62\% with the FS preconditioner; while using 16
processes instead of 8, the reduction is 9.50\% with the AS preconditioner and 8.60\% with the FS preconditioner. 
\new{In general, neither preconditioner shows good scalability properties.
Again the reason for the ``non-perfect'' scalability seems to be the additive Schwarz domain decomposition used between processes.
However, now the negative effect appears to be stronger with respect to the 2D case, since an  
increase in the number of processes causes an consistent increase in the number of linear iterations as well. A similar effect 
is pointed out in \cite{crosetto2011parallel}, where, again, a ``non-perfect'' scalability was noticed.}

Comparing the times required by the two solvers, the times obtained with the FS strategy are considerably smaller 
than those obtained with the AS strategy. For the brain aneurysm case, using the FS preconditioner, the time is reduced by
57.64\% with 4 processes, 57.01\% with 8 processes, and 53.46\% with 16 processes. 
Hence, averaging these three percentage values, there is an overall time reduction of 56.04\%. 
Focusing on the aortic aneurysm case, the FS preconditioner is approximately 1.5 times faster than the AS one. 
Using the FS strategy, the time is reduced by 56.93\% with 4 processes, 53.80\% with 8 processes, and 53.34\% with 16 processes. 
The overall time reduction is 54.69\%. 

\begin{figure}[!tb]
  \begin{center}
     \includegraphics[scale=0.8]{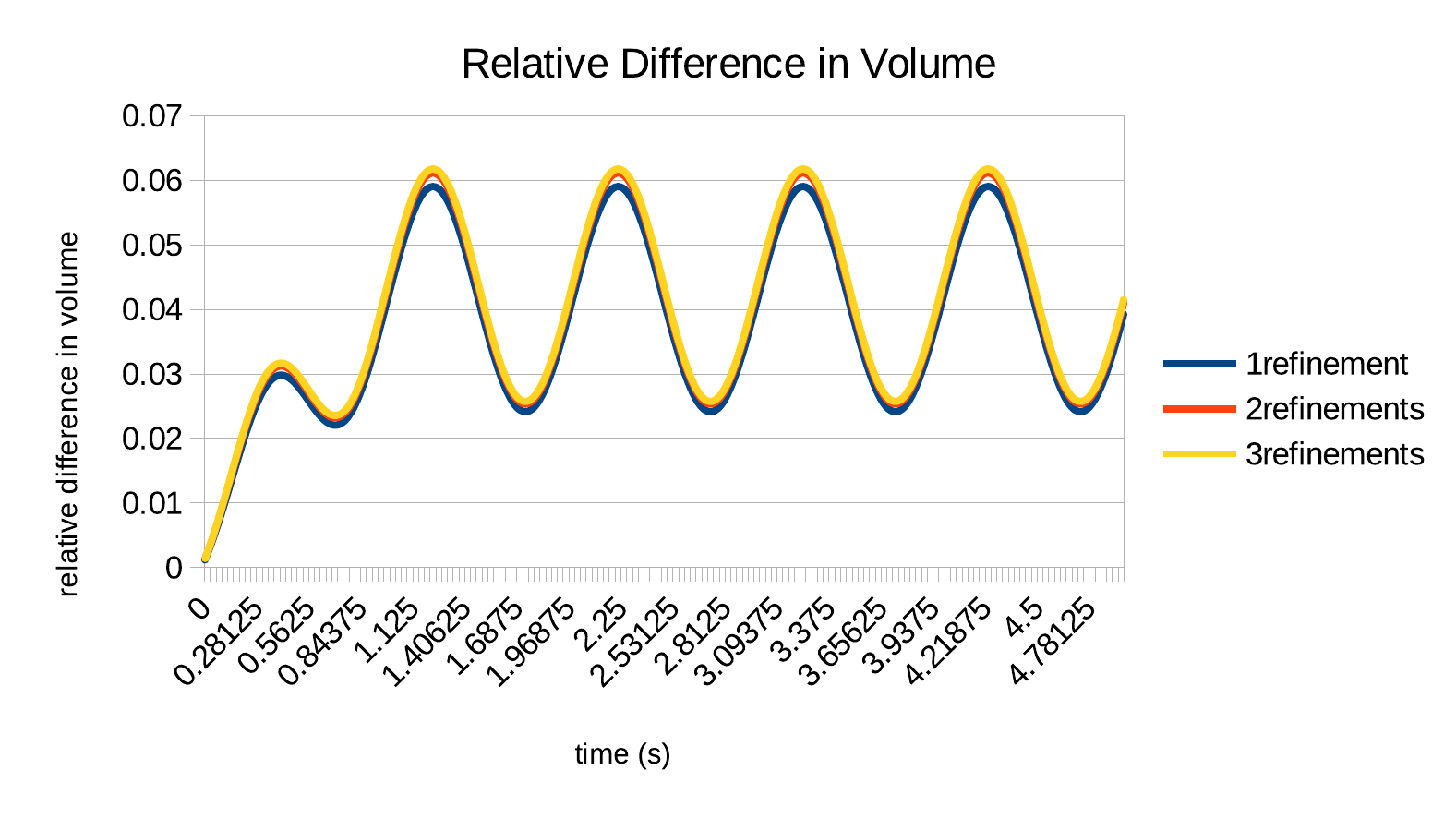}
    \caption{Relative difference in volume of the aneurysm cavity.}\label{volume}
  \end{center}
\end{figure}

%
\new{Once again, we show that the proposed FSI modeling  and solvers are able to correctly capture some physical quantities of
interest. Results obtained with the AS preconditioner are identical, so they are not shown. Focusing on the brain aneurysm, 
Figure \ref{volume} shows the relative difference in volume of the aneurysm cavity. 
This quantity of interest is computed for three different levels of mesh refinements, namely 1, 2 and 3 levels.
Due to the pulsatile nature of the blood stream and blood pressure, we expect the aneurysm cavity to periodically inflate, 
and this behavior is clearly visible in the graph. The discrepancies in the solutions for the three levels of refinements 
are negligible, therefore the solution is consistently represented. 
For this test, a time-dependent normal stress boundary condition of the form 
$$ \boldsymbol (\bm{\sigma}(\bm{u},p^f) \cdot \mathbf n) \cdot \mathbf n= 2500 + 1000 \sin(2\pi t), \quad \boldsymbol (\bm{\sigma}(\bm{u},p^f) 
\cdot \mathbf n) \cdot \boldsymbol\tau= 0$$ 
has been imposed on the outlet, to better simulate the more physiological pulsatile nature of blood pressure.} 

\begin {table}[!tb]
\setlength\tabcolsep{5.5pt} 
\caption{3D Brain Aneurysm: solver tests - 3 refinements}
\begin{center}
	\begin{tabular}{|c|c|c|c|c|c|c|c|c|} \hline	
	\multirow{2}{*} {}   & \multicolumn{8}{|c|}{3D Brain Aneurysm} \\ \cline {2-9}
	& \multicolumn{4}{|c|} {Additive Schwarz } & \multicolumn{4}{|c|} {Field-Split } \\ \hline		      
	 procs  & $s^{max}$ &  $N$  & $\rho$ & Time        & $s^{max}$ &  $N$  & $\rho$ & Time\\ \hline 
	     4  &  3.78     & 20.87 &  0.40  & 11,662.29s  &  3.78     & 20.88 &  0.40  & 4,939.75s \\ \hline 
	     8  &  3.78     & 29.08 &  0.50  & 7,049.10s   &  3.78     & 29.14 &  0.51  & 3,280.01s \\ \hline 
	    16  &  3.81     & 36.19 &  0.67  & 4,957.86s   &  3.84     & 37.00 &  0.69  & 2,307.21s  \\ \hline 
	\end{tabular}
\end{center}
\label{turek3D_ksp&time}
\end{table}

\begin {table}[!tb]
\setlength\tabcolsep{5.5pt} 
\caption{3D Aortic Aneurysm: solver tests - 3 refinements}
\begin{center}
	\begin{tabular}{|c|c|c|c|c|c|c|c|c|} \hline	
	\multirow{2}{*} {}   & \multicolumn{8}{|c|}{3D Aortic Aneurysm} \\ \cline {2-9}
	& \multicolumn{4}{|c|} {Additive Schwarz } & \multicolumn{4}{|c|} {Field-Split } \\ \hline		      
	 procs  & $s^{max}$ &  $N$  & $\rho$ & Time       & $s^{max}$ &  $N$  & $\rho$ & Time\\ \hline 
	     4  &  3.09     & 28.40 &  0.55  & 8,848.78s  &  3.09     & 28.41 &  0.55  & 3,810.90s \\ \hline 
	     8  &  4.09     & 33.21 &  0.70  & 7042.69s   &  4.13     & 33.64 &  0.71  & 3253.45s \\ \hline 
	    16  &  4.66     & 39.11 &  0.81  & 6,373.04s  &  4.78     & 39.46 &  0.82  & 2,973.71s \\ \hline 
	\end{tabular}
\end{center}
\label{AAA3D_ksp&time}
\end{table}

\section{Conclusions}\label{conclusions}
In this paper, a GMRES algorithm preconditioned by GMG is used to solve linearized FSI systems. 
To ensure numerical stability, a monolithic formulation of the FSI system is adopted, and a 
SUPG stabilization is added to the momentum equation. 
The novelty of the work consists of a field-split preconditioner for the multigrid level sub-solvers. 
The block structure of this FS preconditioner derives from using the physical variables as a splitting strategy. The FS preconditioner is a combination of physics-based and domain decomposition preconditioners, since its diagonal blocks are further preconditioned using a locally multiplicative AS strategy. 
The performance of the proposed field-split preconditioner is tested on biomedical FSI applications and compared with those of a pure AS preconditioner. Both 2D and 3D simulations are carried out on aneurysm and venous valve geometries. 
\new{The numerical tests show a weak mesh dependence of the FS and AS preconditioners, as expected with complex biomedical geometries.
The computational times show optimal complexity with respect to the number of unknowns.
Moreover, both solvers are found robust with respect to the time step size and have similar convergence properties.  
The major difference between the FS and the AS preconditioners is that FS is significantly faster than AS, especially in 3D. 
In the examples considered, the FS preconditioner is approximately 20\% faster than the AS in 2D simulations and 50\% faster in 3D. 
For the 2D cases, a comparison is also made against a direct LU solver, showing an even better speed-up.  
In terms of process scalability, the performance of both preconditioners is satisfactory in 2D, but exhibits saturation in 3D.}

Future work consists of investigating the performance of the proposed field-split technique as a preconditioner for the GMRES algorithm. 
Specifically, referring to Figure \ref{FS}, we are interested in exploring how the solver's performance would be affected interchanging the GMG and the FS boxes. 
This interchange would produce a Krylov subspace solver preconditioned with FS, where the diagonal blocks are further preconditioned with a GMG algorithm. 
The smoothing procedure would then be of domain decomposition type (AS). 

\section*{Data Accessibility}
The data associated with this paper has been obtained using FEMuS, an open-source finite element C++ library. 
FEMuS can be downloaded from GitHub; please see the link in reference \cite{femus-web-page}. \\

\section*{Acknowledgments}
The research of the second author was partially supported by the NSF grant DMS-1912902.

\bibliographystyle{plain}
\bibliography{biblio}

\end{document}